\documentclass[11pt]{article}

\usepackage{latexsym}
\usepackage{hyperref}
\usepackage[all]{xy}

\usepackage{amsmath}
\usepackage{xypic}

\usepackage{amsmath}
\usepackage{amssymb}
\usepackage{amsfonts}
\usepackage{amsthm}
\usepackage{mathrsfs}
\usepackage{stmaryrd}
\usepackage[utf8]{inputenc}
\usepackage{xcolor}
\usepackage{newunicodechar}
\usepackage{xspace}

\usepackage{enumitem}
\setlist{nosep} 

\usepackage[lmargin=1.5in,rmargin=1.5in,tmargin=1.5in,bmargin=1.5in]{geometry}

\usepackage{multicol}
\DeclareUnicodeCharacter{00A0}{ }

\newunicodechar{α}{\ensuremath{\alpha}}
\newunicodechar{β}{\ensuremath{\beta}}
\newunicodechar{χ}{\ensuremath{\chi}}
\newunicodechar{δ}{\ensuremath{\delta}}
\newunicodechar{ε}{\ensuremath{\varepsilon}}
\newunicodechar{Δ}{\ensuremath{\Delta}}
\newunicodechar{η}{\ensuremath{\eta}}
\newunicodechar{γ}{\ensuremath{\gamma}}
\newunicodechar{Γ}{\ensuremath{\Gamma}}
\newunicodechar{ι}{\ensuremath{\iota}}
\newunicodechar{κ}{\ensuremath{\kappa}}
\newunicodechar{λ}{\ensuremath{\lambda}}
\newunicodechar{Λ}{\ensuremath{\Lambda}}
\newunicodechar{μ}{\ensuremath{\mu}}
\newunicodechar{ω}{\ensuremath{\omega}}
\newunicodechar{Ω}{\ensuremath{\Omega}}
\newunicodechar{π}{\ensuremath{\pi}}
\newunicodechar{φ}{\ensuremath{\phi}}
\newunicodechar{Φ}{\ensuremath{\Phi}}
\newunicodechar{ψ}{\ensuremath{\psi}}
\newunicodechar{Ψ}{\ensuremath{\Psi}}
\newunicodechar{ρ}{\ensuremath{\rho}}
\newunicodechar{σ}{\ensuremath{\sigma}}
\newunicodechar{Σ}{\ensuremath{\Sigma}}
\newunicodechar{τ}{\ensuremath{\tau}}
\newunicodechar{θ}{\ensuremath{\theta}}
\newunicodechar{Θ}{\ensuremath{\Theta}}
\newunicodechar{ζ}{\ensuremath{\zeta}}
\newunicodechar{ν}{\ensuremath{\nu}}

\newunicodechar{𝔸}{\ensuremath{\mathbb{A}}}
\newunicodechar{ℂ}{\ensuremath{\mathbb{C}}}
\newunicodechar{ℕ}{\ensuremath{\mathbb{N}}}
\newunicodechar{ℙ}{\ensuremath{\mathbb{P}}}
\newunicodechar{ℚ}{\ensuremath{\mathbb{Q}}}
\newunicodechar{ℤ}{\ensuremath{\mathbb{Z}}}

\newunicodechar{𝒪}{\ensuremath{\mathcal{O}}}

\newunicodechar{∞}{\ensuremath{\infty}}
\newunicodechar{→}{\ensuremath{\to}}
\newunicodechar{←}{\ensuremath{\leftarrow}}
\newunicodechar{↦}{\ensuremath{\mapsto}}
\newunicodechar{≅}{\ensuremath{\cong}}
\newunicodechar{×}{\ensuremath{\times}}
\newunicodechar{∪}{\ensuremath{\cup}}
\newunicodechar{∩}{\ensuremath{\cap}}
\newunicodechar{∏}{\ensuremath{\prod}}
\newunicodechar{⊔}{\ensuremath{\sqcup}}
\newunicodechar{⊓}{\ensuremath{\sqcap}}
\newunicodechar{⊇}{\ensuremath{\supseteq}}
\newunicodechar{⊃}{\ensuremath{\supset}}
\newunicodechar{⊆}{\ensuremath{\subseteq}}
\newunicodechar{⊂}{\ensuremath{\subset}}
\newunicodechar{≥}{\ensuremath{\geq}}
\newunicodechar{≤}{\ensuremath{\leq}}
\newunicodechar{∈}{\ensuremath{\in}}
\newunicodechar{◦}{\ensuremath{\circ}}
\newunicodechar{⊕}{\ensuremath{\oplus}}
\newunicodechar{○}{\ensuremath{\circ}}

\DeclareFontFamily{U}{wncy}{}
\DeclareFontShape{U}{wncy}{m}{n}{<->wncyr10}{}
\DeclareSymbolFont{mcy}{U}{wncy}{m}{n}
\DeclareMathSymbol{\Sh}{\mathord}{mcy}{"58} 

\newtheorem{theo}{Theorem}[section]
\newtheorem{prop}[theo]{Proposition}
\newtheorem{lem}[theo]{Lemma}
\newtheorem{rmk}[theo]{Remark}
\newtheorem{defi}[theo]{Definition}
\newtheorem{coro}[theo]{Corollary}
\newtheorem{lemm}[theo]{Lemma}

\theoremstyle{definition}

\newtheorem{exam}[theo]{Example}
\newtheorem{rema}[theo]{Remark}

\newtheorem{claim}[theo]{Claim}

\DeclareMathOperator{\uSpec}{\underline{Spec}}
\newcommand{\cA}{\mathcal{A}}
\renewcommand{\AA}{\mathbb{A}}

\newcommand{\cC}{\mathcal{C}}
\newcommand{\cI}{\mathcal{I}}
\newcommand{\cJ}{\mathcal{J}}

\renewcommand{\cD}{\mathcal{D}}
\newcommand{\cE}{\mathcal{E}}
\newcommand{\FF}{\mathbb{F}}
\newcommand{\GG}{\mathbb{G}}

\newcommand{\cO}{\mathcal{O}}

\newcommand{\sU}{\mathcal{U}}

\newcommand{\NN}{\mathbb{N}}

\newcommand{\OO}{\mathcal{O}}

\newcommand{\QQ}{\mathbb{Q}}

\newcommand{\cS}{\mathcal{S}}

\newcommand{\cU}{\mathcal{U}}
\newcommand{\cV}{\mathcal{V}}

\newcommand{\procdh}{\mathrm{procdh}}

\newcommand{\ARingfp}{\mathsf{aRing}^{\mathrm{ft}}}
\newcommand{\ARingffp}{\mathsf{aRing}^{\mathrm{fft}}}

\newcommand{\pARing}{\mathsf{paRing}}
\newcommand{\ARing}{\mathsf{aRing}}

\newcommand{\HN}{\operatorname{HC}^-}
\newcommand{\TC}{\operatorname{TC}}

\DeclareMathOperator{\Proj}{Proj}

\DeclareMathOperator{\Ab}{Ab}
\DeclareMathOperator{\Spc}{\cS}

\DeclareMathOperator{\Sp}{Sp}

\DeclareMathOperator{\Spec}{Spec}
\DeclareMathOperator{\Spf}{Spf}

\DeclareMathOperator{\Ker}{Ker}

\DeclareMathOperator{\id}{id}

\newcommand{\Spt}{\mathrm{Spt}}

\DeclareMathOperator{\Hom}{Hom}

\DeclareMathOperator{\fib}{fib}

\newcommand{\ZZ}{\mathbb{Z}}

\newcommand{\fM}{\frak{M}}

\newcommand{\fS}{\frak{S}}

\newcommand{\fU}{\frak{U}}
\newcommand{\fV}{\frak{V}}
\newcommand{\fW}{\frak{W}}
\newcommand{\fX}{\frak{X}}
\newcommand{\fY}{\frak{Y}}
\newcommand{\fZ}{\frak{Z}}

\newcommand{\fm}{\frak{m}}
\newcommand{\cont}{\textrm{cont}}

\newcommand{\Set}{\mathsf{Set}}

\newcommand{\PSh}{\mathsf{PSh}}

\newcommand{\Shv}{\mathsf{Shv}}
\newcommand{\Shvzen}{\mathsf{Shv}_{\procdh}}
\newcommand{\Shvzenc}{\mathsf{Shv}^{\cont}_{\procdh}}

\newcommand{\Sm}{\mathsf{Sm}}

\newcommand{\Sch}{\mathsf{Sch}}

\newcommand{\Schfp}{\mathsf{Sch}^{\mathrm{fp}}}
\newcommand{\SchfpS}{\mathsf{Sch}_S^{\mathrm{fp}}}
\newcommand{\fSch}{\mathsf{fSch}}
\newcommand{\fSchad}{\mathsf{fSch}^{ad}}

\newcommand{\fSchfp}{\mathsf{fSch}^{\mathrm{fp}}}

\newcommand{\fSchS}{\mathsf{fSch}_S}
\newcommand{\fSchfS}{\mathsf{fSch}_\fS}
\newcommand{\AffSchS}{\mathsf{AffSch}_S}

\newcommand{\fSchadfS}{\fSchfS^{\mathrm{ad}}}

\newcommand{\fSchfpfS}{\fSchfS^{\mathrm{ft}}}

\newcommand{\fSchffp}{\fSch^{\mathrm{fft}}}
\newcommand{\fSchffpS}{\fSchS^{\mathrm{fft}}}
\newcommand{\fSchffpfS}{\fSchfS^{\mathrm{fft}}}

\newcommand{\AffSchffpS}{\AffSchS^{\mathrm{fft}}}

\def\fc#1#2{\widehat{#1}_{|#2}}
\newcommand{\aRing}{\mathsf{aRing}}

\newcommand{\RingS}{\mathsf{Ring}_S}

\newcommand{\SchS}{\mathsf{Sch}_S}

\newcommand{\RZ}{\mathsf{RZ}}
\newcommand{\Ind}{\mathrm{Ind}}

\newcommand{\Zar}{\mathsf{Zar}}
\newcommand{\Nis}{\mathsf{Nis}}

\newcommand{\cdh}{\mathsf{cdh}}

\newcommand{\Arr}{\mathsf{Arr}}

\newcommand{\LO}[1]{L\Omega_{#1}}
\newcommand{\hLO}[1]{\widehat{L\Omega}_{#1}}
\newcommand{\hLOH}[2]{\widehat{L\Omega}^{\geq #1}_{#2}}
\newcommand{\Znsyn}{{\ZZ_p(n)^{\operatorname{syn}}}}

\newcommand{\Mod}{\mathsf{Mod}}

\newcommand{\ol}{\overline}

\renewcommand{\FF}{\mathbb{F}}
\renewcommand{\GG}{\mathbb{G}}

\newcommand{\qfor}{\quad\text{for }}

\newcommand{\aff}{\mathrm{aff}}

\newcommand{\qcqs}{\mathrm{qcqs}}

\newcommand{\pft}{\mathrm{pft}}

\newcommand{\eq}[2]{\begin{equation}\label{#1}#2 \end{equation}}

\def\SchftSb{\Sch^{\mathrm{ft}}_{S/I^\infty}}
\def\Schft{\Sch^{\mathrm{ft}}}

\def\sI{\mathcal{I}}

\def\Nil{\mathrm{Nil}}

\newcommand{\Ring}{\mathsf{Ring}}

\def\catprojlim#1{\underset{#1}{``\varprojlim"}}
\def\catinjlim#1{\underset{#1}{``\varinjlim"}}

\def\qaq{\;\text{ and }\;}
\def\rmapo#1{\overset{#1}{\longrightarrow}}

\def\Rinf{R_\infty}

\def\Bl{\mathrm{Bl}}
\def\BlAt{\Bl_{0}(\AA^2)}

\def\BlAtZ{\Bl_0(\AA^2_\ZZ)}
\def\hBlAtS{\widehat{\Bl}_{0}(\AA^2_S)}

\def\Fib{\mathrm{Fib}}
\def\Bl{\mathrm{Bl}}

\def\qaq{\;\text{ and }\;}

\def\qwith{\;\text{with }}

\def\fS{\mathfrak{S}}

\def\Rl{R_\lambda}
\def\Rlb{\overline{R}_\lambda}
\def\Rab{\overline{R}_\alpha}

\def\fc#1#2{\widehat{#1}_{|#2}}
\def\ac#1#2{{#1}^{\mathrm{ad}}_{#2}}

\def\Rl{R_\lambda}

\def\kn{\kappa(\fn)}
\def\Rn{{\Rinf}_\fn}

\def\tU{\widetilde{U}}
\def\tY{\widetilde{Y}}

\def\tfU{\widetilde{\fU}}

\def\tfY{\widetilde{\fY}}

\def\hB{\hat{B}}
\def\hC{\hat{C}}

\def\fn{\mathfrak{n}}

\def\Rg{R_\gamma}
\def\fg{f_\gamma}

\def\hotimes{\hat{\otimes}}

\def\hAS#1{\widehat{\AA}^{#1}_S}
\def\hASo#1{\widehat{\AA}^{#1}_{S|0}}
\def\AS#1{\AA^{#1}_S}
\def\AZ#1{\AA^{#1}_\ZZ}
\def\cAx{\widehat{\cA}_x}

\def\tNil{\mathrm{tNil}}
\def\LAP{\mathrm{LA}^{\pft}}

\def\isom{\overset{\simeq}{\longrightarrow}}

\def\alphab{\overline{\alpha}}
\def\fb{\overline{f}}

\def\fl{f_\lambda}
\def\fbl{\fb_\lambda}
\def\abl{\alphab_\lambda}
\def\Rl{R_\lambda}

\def\fml{\fm_\lambda}
\def\hV{\hat{V}}
\def\hC{\hat{C}}
\def\phib{\overline{\phi}}
\def\tpsi{\tilde{\psi}}
\def\phib{\overline{\phi}}
\def\psib{\overline{\psi}}
\def\tpsi{\tilde{\psi}}
\def\Laml{\Lambda_{\lambda/}}
\def\Ra{R_\alpha}

\def\ModX{\Mod(\fX)}

\def\fSchffpSX{{\fSchffpS}_{/\fX}}
\def\fSchffpX{\fSchffp_{\fX}}
\def\Fc{F^{\cont}}
\def\Kc{K^{\cont}}

\def\Af{A_{\{f\}}}

\newcounter{spec}
{\end{list}}%

\begin{document}
\title{A pro-cdh topology on formal schemes}
\author{Shane Kelly and Shuji Saito}
\maketitle

\begin{abstract}
We introduce a pro-cdh topology on formal schemes and prove 
that the $\infty$-topos of pro-cdh sheaves of spaces has an optimal bound of homotopy 
dimension. This remedies a defect for a pro-cdh topology on schemes introduced in \cite{KS23}.
As an application, we give a topos-theoretic interpretation of Weibel's vanishing of negative K-theory and motivic cohomology of Elmanto and Morrow.
\end{abstract}
\tableofcontents

\section*{Introduction}

This article is a sequel to \cite{KS23}, where a pro-cdh topology for schemes is introduced.  The main motivation was to provide a Grothendieck topology encoding the pro-excision property for abstract blowup squares enjoyed by some important invariants such as algebraic $K$-theory and cotangent complexes and to employ topos theoretic techniques to study those invariants. 
It is defined by modifying the cdh topology in such a way that it sees nilpotents and make it possible to define motivic cohomolgy that sees nilpotents.
Here, we recall its definition.
Let $S$ be a qcqs scheme and $\Schfp_S$ be the category of schemes of finite presentation over $S$.

\begin{defi}(\cite[Def. 1.1]{KS23})\label{def-intro;procdhSch}
The \emph{pro-cdh topology} on $\Schfp_S$ is generated by the following coverings.
\begin{enumerate}
 \item Nisnevich coverings.
 \item Proabstract blowup squares: families of the form
 \[ \{ Z_n \to X\}_{n \in \NN} \sqcup \{Y \to X\} \]
 where $Y \to X$ is a proper morphism of finite presentation which is an isomorphism outside of a closed subscheme $Z_0 \subseteq X$ of finite presentation, and $Z_n = \uSpec \OO_X / \sI_Z^{n+1}$ is the $n$th infinitesimal thickening of $Z_0$.
\end{enumerate}
Let $\Shvzen(\Schfp_S)$ (resp. $\Shvzen(\Schfp_S,\Spc)$) be the topos of pro-cdh sheaves of sets (resp. the $\infty$-topos of pro-cdh sheaves of spaces) on $\Schfp_S$. 
\end{defi}

The pro-cdh topology enjoys a number of pleasant properties.
Here, we recall the following.

\begin{theo}(\cite[Th.1.3]{KS23})\label{HomotopyDimensionSch}
If $S$ has finite valuative dimension $d$ with Noetherian underlying topological space, $\Shv_{\procdh}(\SchfpS, \cS)$ has homotopy dimension $\leq 2d$.
\end{theo}

The above theorem implies that $\Shv_{\procdh}(\Schfp_S, \cS)$ is hypercomplete if $S$ has finite valuative dimension $d$ with Noetherian underlying topological space (\cite[Cor.1.4]{KS23}).
Combined with the conservativity of fibre functors of $\Shv_{\procdh}(\Schfp_S)$ and a characterization of those fibre functors (\cite[Th.1.5 and Pr.1.7]{KS23}), it also provides a topos-theoretic interpretation of the Bass construction  (\cite[Th.1.8]{KS23}): For any Noetherian scheme $X$ of finite Krull dimension $d$, there exists a natural equivalence
\[ (a_{\procdh} τ_{\geq 0}K)(X) \simeq K(X).\]
Here $K(X)$ is the non-connective algebraic $K$-theory of $X$ and $τ_{\geq 0}K(X)$ the connective $K$-theory.

However, Theorem \ref{HomotopyDimensionSch} has a drawback: Due to the appearance of $2d$ instead of $d$\footnote{There exists a Noetherian scheme of Krull dimension one with pro-cdh homotopy dimension two, see \cite[Ex.7.17]{KS23}.},
it is not possible to give a topos-theoretic interpretation of Weibel's vanishing of negative $K$-theory, which is the celebrated theorem \cite[Th.B]{KST-Weibel}:
For any Noetherian scheme $X$ of finite Krull dimension $d$, we have   
\eq{KtheoryWeibel-intro}{K_i(X)=0\qfor i<-d.}

A main purpose of this article is to remedy the drawback by considering a pro-cdh topology on formal schemes. We now introduce our basic setup. 
Let $\ARing$ be the category of adic rings with continuous ring homomorphisms (see Definition \ref{def;adicring}).
Let $\fSch$ (resp. $\fSch_S$) denote the category of formal schemes (resp. formal schemes over a fixed base formal scheme $S$). 
In this paper, we only treat locally Noetherian formal schemes $\fX$ meaning that any point of $\fX$ admits an affine formal open neighborhood $\Spf(A)$ with $A\in \ARing$ Noetherian.

Let $\fX$ be a locally Noetherian formal scheme and $f:\fY\to \fX$ be a morphism in $\fSch$.
We say that $f$ is \emph{of finite type} (resp. \emph{formally of finite type}) if  it is quasi-compact and if for any point $y\in \fY$, there exists an affine formal open neighborhood $U=\Spf(A)$ of $f(y)$ in $\fX$ and an affine formal open neighborhood $V=\Spf(B)$ of $y$ in $\fY$ such that $f(V)\subset U$ and $B$ is \emph{topologically of finite type over $A$} (resp. \emph{formally of finite type over $A$}). Here, a map $A\to B$ in $\ARing$ is called \emph{topologically of finite type} (resp. \emph{formally of finite type}) if $B$ is a quotient of a toplogical $A$-algebra of the form $ A\{z_1,\dots,z_r\}$ (resp. $ A\{z_1,\dots,z_r\}[[w_1,\dots,w_s]]$), where $ A\{z_1,\dots,z_r\}$ is the convergent power series ring over $A$, see Definition \ref{def;ftAring}.

For a locally Noetherian formal $S$, we let $\fSchffpS$ denote the full subcategory of $\fSch_S$ of formal schemes formally of finite type over $S$.

\begin{defi} \label{defi-intro:procdh}
The \emph{pro-cdh topology} on $\fSchffpS$ is generated by the following coverings.
\begin{enumerate} 
 \item Nisnevich coverings, Definition \ref{def;fXnis}.
\item Formal abstract blowup covering: families of the form
 \[ \{\fc \fX Z \sqcup \fY \to \fX\} \]
 where $Z$ is a closed subscheme of $\fX$ and $\fc \fX Z$ is the formal completion of $\fX$ along $Z$, see Definition \ref{def;formalCompletion},
and $f:\fY \to \fX$ is proper
and of finite type and an isomorphism over the open complement of $Z$ in $\fX$.
Moreover, $f$ is \emph{locally algebraizable}, see Definition \ref{defi:locallyalgebriaable}.
\end{enumerate}
We let $\Shvzen(\fSchffpS)$ (resp. $\Shvzen(\fSchffp_S,\Spc)$) denote the topos of pro-cdh sheaves of sets (resp. the $\infty$-topos of pro-cdh sheaves of spaces) on $\fSchffp_S$. 
\end{defi}

We will prove the following result on the bound of the homotopy dimension of $\Shvzen(\fSchffp_S,\Spc)$, which improves the bound from Theorem \ref{HomotopyDimensionSch}.

\begin{theo}[Theorem \ref{thm2;Hypercomplete}]\label{thm2Intro;Hypercomplete}
If the underlying topological space of $S$ has finite Krull dimension $d$, 
$\Shvzen(\fSchffpS,\cS)$ has homotopy dimension$\leq d$.
\end{theo}

Now, we explain how the pro-cdh topology on schemes is related to that on formal schemes, see Remark \ref{rem;defi:procdh} for more detail.
Let $S$ be a locally Noetherian scheme viewed as a discrete formal scheme.
Let $\Schft_S$ be the category of schemes of finite type over $S$.
By definition, $\Schft_S$ is viewed as a full subcategory of $\fSchffpS$.
Then, the right Kan extension and the restriction along $\Schft_S\to \fSchffpS$ induce an equivalence
\eq{Shvprocdhcont-intro}{\Shvzen(\Schft_S,\Spc) \simeq \Shvzenc(\fSchffpS,\Spc),}
where the left hand side is from Definition \ref{def-intro;procdhSch} and the right hand side is the full subcategory of $\Shvzen(\fSchffpS,\Spc)$ of those sheaves satisfying the condition
\eq{Cont-condition}{F(\fX) = \varprojlim_n F(\fX_n)\qfor \fX\in \fSchffpS,}
where $\fX_n=\uSpec(\cO_\fX/\cI^n)$ for an ideal $\cI$ of definition of $\fX$  (\cite[10.3.1]{EGAI}).
The condition \eqref{Cont-condition} is equivalent to the sheaf condition for the topology on $\fSchffpS$ whose covering families are those of the form: $\{ \fX_n \to \fX \}_{n \in \NN}$ for $\fX\in\fSchffpS$ .
Hence, $\Shvzen(\Schft_S,\Spc)$ is a topological localization of $\Shvzen(\fSchffpS,\Spc)$.
\medbreak

In order to apply Theorem \ref{thm2Intro;Hypercomplete} to a topos-theoretic interpretation of Weibel's vanishing, we need to characterize the fibre functors of $\Shvzen(\fSchffpS)$, i.e.
functors $ \Shvzen(\fSchffpS)\to \Set$ which preserves colimits and finite limits. 
Let $S$ be a Notherian adic ring and let $\ARingffp_S$ be the category of adic rings formally of finite type over $S$. Write $\fSchffpS$ for $\fSchffp_{\Spf(S)}$.

\begin{theo}[Theorem \ref{prop:procdhLocal}]\label{prop-intro:procdhLocal}
There exists a bijection between the fibre functors of $\Shvzen(\fSchffpS)$ and the pairs $(R,\fn)$ of a ring $R$ and an ideal $\fn\subset R$ satisfying the following conditions:
\begin{enumerate}
\item[(i)]
There exists an ind-object $R_\bullet=\catinjlim {\lambda\in \Lambda}  \Rl $ in $\ARingffp_S$ such that 
\[R=\varinjlim_\lambda \Rl \qaq \fn=\varinjlim_{\lambda}\fn_\lambda,\]
where $\fn_\lambda\subset R_\lambda$ is the ideal of topologically nilpotent elements and these colimits take place in the category of rings.
\item[(ii)]
$R/\fn$ is a henselian valuation ring and
\[R\simeq R/\fn\times_{\kn} R_\fn,\]
where $\kn$ is the residue field of $\fn$ and $R_\fn$ is the localization of $R$ at $\fn$.
\end{enumerate}
\end{theo}

The fibre functor associated to a pair $(R,\fn)$ is given by
\eq{phiR-intro}{ \phi_{R_\bullet}: \Shvzen(\fSchffp_S) \to \Set\;;\; F \to 
\varinjlim_{\lambda\in \Lambda} F(\Spf(R_\lambda))}
for a chosen ind-object $R_\bullet$ as in (i) and it is independent of choices, Proposition \ref{prop;fibrePsi}.
\medbreak

Finally, using Theorems \ref{thm2Intro;Hypercomplete} and Theorem~\ref{prop-intro:procdhLocal}, we will give a topos-theoretic new proof of the following axiomatization of Weibel's vanishing \eqref{KtheoryWeibel-intro}, which was proved by Elmanto-Morrow \cite[Prop.8.10]{EM23} using an argument of Kerz-Strunk-Tamme for their proof of Weibel's conjecture \cite[Th.B]{KST-Weibel} and also an idea 
of Cortiñas–Haesemeyer(–Schlichting)–Weibel for their proofs of Weibel and Vost's conjecture for finite type schemes over fields of characteristic zero 
\cite{CHSW08}, \cite{CHW08}.

We set $\cC = \Spt$ or $D(\ZZ)$. For an integer $N$, $\cC_{\geq N}$ denotes the full subcategory of objects of $\cC$ supported in homological degree $i\geq N$.  
Let $\Sch^{\qcqs}$ be the category of qcqs schemes.

\begin{theo}[Theorem \ref{theo:procdhWebelVanishing}] \label{theo-intro:procdhWebelVanishing}
For $F \in \PSh(\Sch^{\qcqs}, \cC)$ and integer $N$, consider the conditions.%
\begin{enumerate}
 \item[(Desc)] For every Noetherian scheme $X$, the restriction $F$ to the category of schemes of finite type over $X$ is a pro-cdh sheaf in the sense of Definition \ref{def-intro;procdhSch}.
 \item[(Fin)] $F$ is finitary, in the sense that it preserves filtered colimits of rings.
 \item[(Rig)$_N$] 
$\fib(F(A)\to F(A/I))\in \cC_{\geq -N}$ for every Noetherian ring $A$ and nilpotent ideal $I\subset A$.
\item[(Val)$_N$] $F(R)\in \cC_{\geq -N}$ for every henselian valuation ring $R$.
\end{enumerate}
If $F$ satisfies (Desc), (Fin), (Rig)$_N$, and (Val)$_N$,
then for every Noetherian scheme $X$ of finite Krull dimension $d$, we have
\begin{equation*} \label{eq;rmk;apcdhWebel}
\pi_i(F(X))=0 \qfor  i<-N-d.
\end{equation*}
\end{theo}

Thanks to \cite[Th.A]{KST-Weibel} and \cite[Th.1.3(iii)]{KM21},
non-connective $K$-theory satisfies (Desc), (Fin), (Rig)$_0$ and (Val)$_0$.
Thus, \eqref{KtheoryWeibel-intro} is a consequence of Theorem \ref{theo-intro:procdhWebelVanishing}.
In \S\ref{Weibel}, we will produce a number of presheaves that satisfy the conditions of Theorem~\ref{theo-intro:procdhWebelVanishing}, using negative cyclic homology, integral topological cyclic homology and motivic cohomology defined by Elmanto and Morrow.
These vanishing results are not new and known from \cite{EM23}.
Theorem \ref{theo-intro:procdhWebelVanishing} gives a new proof of these results.

\subsection*{Acknowledgements}
The authors are deeply grateful to J. Ayoub for his suggestion on 
a substantial simplification of the proof of Proposition \ref{prop:cocont}.
The authors also thank M. Morrow for helpful comments.
The second author thanks A. Abbes for answering several questions on formal geometry. 


\section{Preliminaries from formal geometry}\label{ClosedPairs}
\def\het#1{{#1}^h}
\def\zet#1{{#1}^{\Zar}}
\def\otz{\otimes^{\Zar}}

We introduce the basic framework where we develop a pro-cdh topos theory.
For a commutative ring $S$, let $\RingS$ be the category of (commutative) rings over $S$ and $\SchS$ be the category of $S$-schemes.

\subsection{Adic rings}\label{ReviewAdicRing}

\def\ac#1#2{#1^\wedge_{#2}}

\begin{defi}\label{def;adicring}
A \emph{preadic ring} is a ring $R$ equipped with the $I$-adic topology for an ideal $I\subset R$. 
By definition, $\{a+I^n|\; a\in R,\; n=1,2,\dots\}$ is a basis for this topology.
An ideal $J\subset R$ is called an \emph{ideal of definition} if there exist $m,n>0$ such that $I^m\subset J^n\subset I$.
An \emph{adic ring} is a preadic ring which admits a finitely generated ideal $I$ of definition and is Hausdorff and complete for the $I$-adic topology.  

Let $\pARing$ be the category of preadic rings with continuous ring homomorphisms and $\ARing\subset \pARing$ be the full subcategory of adic rings.
\end{defi}

\def\fg{\mathrm{fg}}

\begin{rmk}\label{rem;def;topNil} 
An element $a\in A$ is \emph{topologically nilpotent} if for any ideal $I\subset A$ of definition, $a^n\in I$ for some $n>0$. 
For $A\in \pARing$, the radical $\fn$ of any ideal $I\subset A$ of definition is 
the ideal of all topologically nilpotent elements, and $\fn$ is the maximal ideal of definition if $\fn$ is finitely generated, in particular if $A$ is Noetherian.
\end{rmk}

\begin{rmk}\label{rem;aRignNotherian}  
Let $A$ be a Noetherian adic ring with an ideal $I$ of definition.
Then, any ideal $J\subset A$ is closed and  $A/J$ equipped with $I(A/J)$-adic topology is an adic ring, \cite[Prop. 7.4.14, 7.4.16, 7.4.12]{FK}.
\end{rmk}

\begin{rmk}\label{rem0;def;adicring} 
For $A,B\in \pARing$, a ring homomorphism $f:A\to B$ is continuous if and only if for every ideal $J\subset B$ of definition, there is an ideal $I\subset A$ of definition such that $f(I)\subset J$. If the ideal $\fn\subset B$ of topologically nilpotent elements is an ideal of definition, then $f$ is continuous if and only if $f(I)\subset \fn$ for 
all  ideals $I\subset A$ of definition.
\end{rmk}

\begin{rmk}\label{rem;def;adicring} 
For $R\in \ARing$, we have a natural isomorphism of topological rings
\[ R\simeq \varprojlim_n R/I^n\qfor R\in \ARing ,\]
where $R/I^n$ is given the discrete topology.
For $A\in \pARing$ with a finitely generated ideal $J$ of definition, 
we have bijections
\eq{eq1;rem;def;adicring}{\Hom_{\pARing}(A,R) \simeq 
\varprojlim_n \varinjlim_m \Hom_{\Ring_S}(A/J^m,R/I^n)
\simeq \Hom_{\ARing}(\fc A J,R) ,}
where $\fc A J=\varprojlim_n A/J^n$ is the completion of $A$ with respect to $J$.
\end{rmk}

\begin{rmk}\label{rem2;def;adicring}(\cite[10.7.2]{EGAI}) 
The category $\ARing$ admits fiber coproducts: For $A\to B$ and $A\to C$ in $\ARing$ with ideals of definition $J\subset B$ and $L\subset C$, their fiber coproduct over $A$ is defined by 
\eq{hotimes}{ B\hotimes_A C =\varprojlim_{n} (B/J^n\otimes_A C/L^n)
=\varprojlim_{n} B\otimes_A C/(J\otimes_A C+ B\otimes_A L)^n.}
Note that the image of $J\otimes_A C + B\otimes_A L$ in $B\hotimes_A C$ generates an ideal of definition. 
\end{rmk}



\begin{defi}\label{def;adicring}
A morphism $A\to B$ in $\pARing$ is called \emph{adic} if for some ideal of definition $I$ of $A$, $IB$ is an ideal of definition of $B$. If the condition holds for some $I\subset A$, it holds for any  ideal of definition of $A$. 
\end{defi}

\begin{lem}\label{lem1;adicring}
For an adic morphism $A\to B$ in $\pARing$ and $R\in \pARing$, the following diagram is cartesian.
\[\xymatrix{
\Hom_{\pARing}(B,R)\ar[r]\ar[d] & \Hom_{\pARing}(A,R)\ar[d] \\
\Hom_{\Ring}(B,R)\ar[r] &\Hom_{\Ring}(A,R)\\}\]
\end{lem}
\begin{proof}
Take $g\in \Hom_{\Ring}(B,R)$ whose image in $\Hom_{\Ring}(A,R)$ comes from $\Hom_{\pARing}(A,R)$.
This means that for any ideal $J\subset R$ of definition, there is an ideal $I\subset A$ of definition such that $g(IB)\subset J$.
Since $IB\subset B$ is an ideal of definition, this means $g\in \Hom_{\pARing}(B,R)$, which completes the proof.
\end{proof}

\begin{defi}\label{def;ftAring}
Let $A$ be a Noetherian adic ring and 
$ A\{z_1,\dots,z_r\}$ be the convergent power series ring over $A$ with variables $z_1,\dots,z_r$. Note $A\{z_1,\dots,z_r\}=\varprojlim_{n>0} A/I^n[z_1,\dots,z_r]$ for an ideal $I\subset A$ of definition.
\begin{itemize}
\item[(1)]
We say that an adic ring $B$ over $A$ is \emph{topologically of finite type} if it is isomorphic to $C= A\{z_1,\dots,z_r\}/J$ with the $IC$-adic topology, where $J$ is an ideal in $ A\{z_1,\dots,z_r\}$.
\item[(2)]
We say that an adic ring $B$ over $A$ is \emph{formally of finite type} if it is a quotient of a toplogical $A$-algebra of the form
$C= A\{z_1,\dots,z_r\}[[w_1,\dots,w_s]]$ with the $IC+(w_1,\dots,w_s)$-adic topology.
\end{itemize}
\end{defi}

\begin{defi}\label{def;catftadicring}
Let $S$ be a Noetherian adic ring. We let $\ARingfp_S$ (resp. $\ARingffp_S$) denote the category of adic rings topologically of finite type (resp. formally of finite type) over $S$ with continuous ring homomorphisms.

\end{defi}

\begin{rema}\label{rema2;adicring}
Let $S$ be a Noetherian adic ring. 
\begin{itemize}
\item[(1)]
Since  $S\{z_1,\dots,z_r\}[[w_1,\dots,w_s]]$ is Noetherian (\cite[0.7.5.4.]{EGAI}),
every $A\in \ARingffp_S$ is a Noetherian adic ring (cf. Remark \ref{rem;aRignNotherian}).
\item[(2)]
By definition, for every $A\in \ARingfp_S$ (resp. $A\in \ARingffp_S$), any quotient $\ol{A}$ with $I\ol{A}$-adic topology is in $ \ARingfp_S$ (resp. $\ARingffp_S$)
\item[(3)]
For every $A\in \ARingffp_S$ with an ideal $I$ of definition and $f\in A$, 
let $\Af$ be the $I$-adic completion of $A[f^{-1}]$.
Then, $\Af\in \ARingfp_S$. Note that $\Af=0$ if $f$ is topologically nilpotent.  
\item[(4)]
Both $\ARingfp_S$ and $\ARingffp_S$ are closed under $\hotimes$.
Indeed, this holds without assuming $S$ Noetherian. 
\end{itemize}
\end{rema}

\begin{rema}\label{rema3;adicring}
For $R\in \ARing_S$ with an ideal $J$ of definition, we have by \eqref{eq1;rem;def;adicring}
\begin{multline*}
 \Hom_{\ARing_S}(S\{z_1,\dots,z_r\}[[w_1,\dots,w_s]],R)\\
=\varprojlim_n\varinjlim_{m,l}  
\Hom_{\Ring_S}(S/I^m[z_1,\dots,z_r,w_1,\dots,w_s]/(w_1,\dots,w_s)^l,R/J^n).\end{multline*}
Hence, we have a bijection
\eq{SxyR}{
\Hom_{\ARing_S}(S\{z_1,\dots,z_r\}[[w_1,\dots,w_s]],R) \simeq 
R^{\times r}\times  \tNil(R)^{\times s}
}
sending $f$ to $(f(z_1),\dots,f(z_r),f(w_1,\dots,f(w_s))$,
where $\tNil(R)\subset R$ is the ideal of topologically nilpotent elements.
For $A\in \ARingffp_S$ of the form
\eq{Affpresentation}{A=S\{z_1,\dots,z_r\}[[w_1,\dots,w_s]]/(f_1,\dots,f_t),}
with $f_i\in S\{z_1,\dots,z_r\}[[w_1,\dots,w_s]]$, \eqref{SxyR} implies 
\eq{SxyRxy}{
\Hom_{\ARing_S}(A,R) \simeq \{\alpha\in R^{\times r}\times  \tNil(R)^{\times s}|\; f_i(\alpha)=0\in R\; (i=1,\dots,t)\}.}
\end{rema}

\begin{defi}\label{def;Phi}
Let $\Ind(\ARing)$ be the category of ind-objects in $\ARing$.
Let
\eq{Psi}{ 
\Psi: \Ind(\ARing)\to\pARing\;;\; \catinjlim {\lambda\in \Lambda} \Rl \mapsto (\Rinf,\fn)}
be a functor defined as follows: The underlying ring $\Rinf$ is the colimit $\varinjlim_\lambda \Rl$ in the category of discrete rings. We equip it with an ideal of definition 
\begin{equation}\label{frakn}
 \fn=\varinjlim_{\lambda}\fn_\lambda,\end{equation}
where $\fn_\lambda\subset R_\lambda$ is the ideal of topologically nilpotent elements,
which is the maximal ideal of definition of $\Rl$, Remark \ref{rem;def;topNil}. 
We point out that $\Rinf$ is not in general separated with respect to $\fn$-adic topology. Let $\iota_\lambda:\Rl \to \Rinf$ be the natural map.

\end{defi}

\begin{lem}\label{lem;Phi}
For $(\Rinf,\fn)$ as above, $\Rinf/\fn$ is reduced so that $\fn=\sqrt{\fn}$. Moreover, 
$\Rinf$ is henselian with respect to $\fn$.
\end{lem}
\begin{proof}
By definition we have $\Rinf/\fn\simeq \varinjlim_\lambda \Rg/\fn_\lambda$ and $\Rg/\fn_\lambda$ are reduced so that $\Rinf/\fn$ is reduced. The last assertion follows from \cite[Cor.3.6]{Greco1}.
\end{proof}

\begin{rema}\label{rem1;lem;Phi}
By Remark \ref{rem0;def;adicring}, Lemma \ref{lem;Phi} implies 
\[\Hom_{\pARing}(A,\Rinf)=\{f\in \Hom_{\Ring}(A,\Rinf)|\; f(I)\subset \fn\}\]
 for $A\in \pARing$ and every ideal $I\subset A$ of definition.
\end{rema}

\def\phiR{\phi_{R_\bullet}}

\begin{rema}\label{rem2;lem;Phi}
For $A\in \ARingffp_S$ and $R_\bullet=\catinjlim {\lambda\in \Lambda} \Rl\in \Ind(\ARing_S)$, put
\eq{phiR}{
\phiR(A) =\varinjlim_{\lambda\in \Lambda}\Hom_{\ARing_S}(A,\Rl) .}
Assume that $A$ has the form  as \eqref{Affpresentation} and consider the maps
\eq{eq;rem2;lem;Phi}{
\phiR(A)  \rmapo{\iota} 
\Hom_{\pARing_S}(A,\Rinf) \to \Rinf^{\times r}\times  \Rinf^{\times s},}
where $(\Rinf,\fn)=\Psi(R_\bullet)$ \eqref{Psi}, $\iota$ is induced by $\iota_\lambda:\Rl\to \Rinf$ and the second map is the evaluation at the images of $(z_1,\dots,z_r,w_1,\dots,w_s)$ in $A$.
By \eqref{SxyRxy}, the composite map is injective and its image is the subset of 
elements $\alpha\in \Rinf^{\times r}\times  \fn^{\times s}$ such that 
there exists $\lambda$ and $\alpha_\lambda\in \Rl^{\times r}\times  \fn_\lambda^{\times s}$ with $\alpha=\iota_\lambda(\alpha_\lambda)$ and $ f_i(\alpha_\lambda)=0\in \Rl$ for $i=1,\dots,t$.
In particular,  $\iota$ is injective.

\end{rema}
\def\alpb{\overline{\alpha}}
\def\fnl{\fn_\lambda}
\def\fnlb{\overline{\fn}_\lambda}

\begin{lemm}\label{lem;reduction}
Take $R_\bullet=\catinjlim {\lambda\in \Lambda} \Rl\in \Ind(\ARing_S)$ with
$\Psi(R_\bullet)=(\Rinf,\fn)$.
For each $\lambda\in \Lambda$, let $\Rlb$ be the image of the natural map 
$\iota_\lambda:\Rl\to \Rinf$ equipped with $\fnlb$-adic topology, where 
$\fnlb=\fnl\Rlb$.
By Remark \ref{rem;aRignNotherian}, it gives an object $\ol{R}_\bullet=\catinjlim {\lambda\in \Lambda} \Rlb\in \Ind(\ARing_S)$ such that $\Psi(\ol{R}_\bullet)=(\Rinf,\fn)$.
Then, the natural map $\phiR(A)\to \phi_{\ol{R}_\bullet}(A)$ is bijective 
for every $A\in \ARingffp_S$.
\end{lemm}
\begin{proof}
We have a commutative diagram
\[\xymatrix{
\varinjlim_{\lambda\in \Lambda}\Hom_{\ARing_S}(A,\Rl) 
\ar[r]\ar[d] & \Rinf^{\times r}\times  \Rinf^{\times s} \\
\varinjlim_{\lambda\in \Lambda}\Hom_{\ARing_S}(A,\Rlb) \ar[ru]\\}\]
where the horizontal and diagonal maps are from \eqref{eq;rem2;lem;Phi} and they are injective. By the description of their images from Remark \ref{rem2;lem;Phi}, it suffices to show the following claim: Take $\lambda\in \Lambda$ and $\alpb\in \Rlb^{\times r}\times  \fnlb^{\times s}$ such that $ f_i(\alpb)=0\in \Rlb$ for $i=1,\dots,t$ and let $\alpha\in \Rl^{\times r}\times  \fnl^{\times s}$ be a lift of $\alpb$. 
Then, there exists a transition map $\psi:\Rl\to R_\mu$ such that 
$f_i(\psi(\alpha))=0\in R_\mu$ for $i=1,\dots,t$.
Indeed, we have $f_i(\alpha)\in J:=\Ker(\Rl\to \Rlb)=\Ker(\Rl\to \Rinf)$. Since $\Rl$ is Noetherian, $J$ is finitely generated so that there is $\psi:\Rl\to R_\mu$ such that $J=\Ker(\psi)$. 
Then, $f_i(\psi(\alpha))=\psi f_i(\alpha)=0\in R_\mu$ as wanted.
This completes the proof. 
\end{proof}

\begin{coro}\label{cor;lem;reduction}
Let $A\in \ARingffp_S$ and $R_\bullet, R'_\bullet \in \Ind(\ARingffp_S)$.
If $\Psi(R_\bullet)=\Psi(R'_\bullet)$, there is a canonical bijection
$\phiR(A)\simeq \phi_{R'_\bullet}(A)$. 
\end{coro}
\begin{proof}
Write $R_\bullet=\catinjlim {\lambda\in \Lambda} \Rl$ and $R'_\bullet=\catinjlim {\gamma\in \Gamma} \Rg'$. By Lemma \ref{lem;reduction}, we may assume that 
$\Rl$ and $\Rg'$ are sub-$S$-algebras of $\Rinf$, where $\Psi(R_\bullet)=(\Rinf,\fn)$.
By Remark \ref{rem2;lem;Phi}, both $\phi_{R_\bullet}(A)$ and $\phi_{R’_\bullet}(A)$ are identified with subgroups of $\Rinf^{\times r} \times \fn^{\times s}$.
Moreover, an element of $\phi_{R_\bullet}(A)$ is represented by 
$a\in R_\lambda^{\times r} \times \fn_\lambda^{\times s}$ for some $\lambda$ satisfying 
$f_i(a)=0$ for $i\in[1,r+s]$.
Since $R_\lambda\subset \Rinf=\varinjlim_\gamma R’_\gamma$, there is 
$\gamma\in \Gamma$ such that 
$a\in {R'_\gamma}^{\times r} \times {\fn'_\gamma}^{\times s}$ so that 
$a\in \phi_{R’_\bullet}(A)$. Thus, we get $\phi_{R_\bullet}(A)=\phi_{R’_\bullet}(A)$.
\end{proof}

\begin{lem}\label{lem;cartesian}
Let the notation be as in Definition \ref{def;Phi} and assume $\Rl$ are Noetherian.
For $A\in \pARing$, the following diagram is cartesian.
\[\xymatrix{\varinjlim_\lambda \Hom_{\pARing}(A,R_\lambda)\ar[r]\ar[d]^{\subset} 
& \Hom_{\pARing}(A,\Rinf)\ar[d]^{\subset} \\
\varinjlim_\lambda \Hom_{\Ring}(A,R_\lambda)\ar[r] &\Hom_{\Ring}(A,\Rinf)\\}\]
\end{lem}  
\begin{proof}
Take $\lambda\in\Lambda$ and $f\in \Hom_{\Ring}(A,R_\lambda)$ such that $g=\iota_\lambda\circ f\in\Hom_{\pARing}(A,\Rinf)$, where $\iota_\lambda:\Rl\to \Rinf$ is the natural map. 
It suffices to find a transition map $t:\Rl\to R_\alpha$ such that $t\circ f$ is in $\Hom_{\pARing}(A,R_\alpha)$. 
We have $g(I)\subset\fn$, Remark \ref{rem1;lem;Phi}.
Since $\fn=\varinjlim_{\mu\in \Lambda} \fn_\mu$ and $I$ is finitely generated, there is 
$s: R_\lambda\to R_\mu$ such that $g(I)\subset \iota_\mu(\fn_\mu)$. 
Noting $g=\iota_\lambda\circ f=\iota_\mu\circ s\circ f$, we get 
\[s\circ f(I) \subset \iota_\mu^{-1}(\iota_\mu(\fn_\mu))=\fn_\mu+K_\mu\qwith K_\mu=\Ker(\iota_\mu).\]
Since $R_\mu$ is Noetherian, $K_\mu$ is finitely generated so that there is $s':R_\mu\to R_\alpha$ such that $K_\mu=\Ker(s')$. Putting $t=s'\circ s$, we get
\[t \circ f(I)\subset s'(\fn_\mu+K_\mu)
=s'(\fn_\mu)\subset \fn_\alpha,\]
which proves the desired assertion.
\end{proof}

The following proposition may be viewed as `relative compactness' for an object of 
$\ARingffp_S$ and will play a key role in the proof of Theorem \ref{prop:procdhLocal}.

\begin{prop}\label{lem;compactness}
Let $S$ be a Noetherian adic ring and take $A\in \ARingffp_S$ with an ideal $I\subset A$ of definition. 
Let $B$ be an $A$-algebra of finite presentation 
and $\hB$ be the $IB$-adic completion $\hB$ of $B$.
For an ind-object $\catinjlim {\lambda\in \Lambda} \Rl$ in $\ARingffp_S$, consider the commutative diagram
\[\xymatrix{
\varinjlim_\lambda \Hom_{\ARing_S}(\hB,R_\lambda)\ar[r]^{\simeq}\ar[d]^{p}
&\varinjlim_\lambda \Hom_{\pARing_S}(B,R_\lambda)\ar[r]^-{\iota}
&\Hom_{\pARing_S}(B,\Rinf)\ar[d]^q\\
\varinjlim_\lambda \Hom_{\ARing_S}(A,R_\lambda)\ar[rr]^{\iota}
&&\Hom_{\pARing_S}(A,\Rinf)\\}\]
where 
$p$ and $q$ are induced by $A\to B\to \hB$ and the bijection comes from \eqref{eq1;rem;def;adicring}.
Then, $f\in \varinjlim_\lambda \Hom_{\ARing}(A,R_\lambda)$ lies in the image of $p$ if and only if $\phi=\iota(f)\in \Hom_{\ARing}(A,\Rinf)$ lies in the image of $q$.
\end{prop}  
\begin{proof}
Choose $\lambda\in \Lambda$ and $f\in \Hom_{\ARing_S}(A,R_\lambda)$ and put
$\phi=\iota(f)\in \Hom_{\pARing_S}(A,\Rinf)$. We view $\Rinf$ as an object of $\pARing_A$ via $\phi$. 
Let $\Laml\subset \Lambda$ be the subcategory of objects under $\lambda$. 
Let $\Rlb$ be as in Lemma \ref{lem;reduction}.
If $\gamma: \Rlb\to \Rab$ is a transition map, we have 
$\overline{\iota}_\alpha\circ \gamma=\overline{\iota}_\lambda$, where $\overline{\iota}_\alpha: \Rab\to \Rinf$ is the natural map.
Since $\overline{\iota}_\alpha$ is injective, $\gamma$ must be unique.
Thus, we can view $\Rab$ for $\alpha\in \Laml$ as an object of $\ARing_A$ using $f$. 
Then, we consider the following diagram, where all squares are commutative:
\[\xymatrix{
\varinjlim_{\alpha\in \Laml} \Hom_{\ARing_A}(\hB,\Rab)\ar[r]^-{\simeq}\ar[d]^{\subset} 
&\varinjlim_{\alpha\in \Laml} \Hom_{\pARing_A}(B,\Rab)\ar[r]^-{\iota}
&\Hom_{\pARing_A}(B,\Rinf)\ar[d]^{\subset} \\
\varinjlim_{\alpha\in \Laml} \Hom_{\ARing_S}(\hB,\Rab)\ar[r]^-{\simeq}
&\varinjlim_{\alpha\in \Laml} \Hom_{\pARing_S}(B,\Rab)\ar[r]^-{\iota}
&\Hom_{\pARing_S}(B,\Rinf)\\
\varinjlim_{\alpha\in \Laml} \Hom_{\ARing_S}(\hB,\Ra)\ar[r]^-{\simeq}\ar[u]^\pi_{\simeq}\ar[d]^p
&\varinjlim_{\alpha\in \Laml} \Hom_{\pARing_S}(B,\Ra)\ar[r]^-{\iota}
&\Hom_{\pARing_S}(B,\Rinf)\ar[u]^{=}\ar[d]^q\\
\varinjlim_{\alpha\in \Laml} \Hom_{\ARing_S}(A,\Ra)\ar[rr]^-{\iota}
&&\Hom_{\pARing_S}(A,\Rinf)\\}\]
where $\pi$ is induced by $\Ra\to \Rab$ and it is a bijection by Lemma \ref{lem;reduction} and the horizontal bijections come from \eqref{eq1;rem;def;adicring}.
Now, assume that there exists $\psi\in \Hom_{\pARing_S}(B,\Rinf)$ such that $\phi=q(\psi)$.
Note $\psi\in \Hom_{\pARing_A}(B,\Rinf)$ by definition.
By Claim \ref{claim;lem;compactness}, it lifts to $\overline{g} \in \varinjlim_{\alpha\in \Laml} \Hom_{\ARing_A}(\hB,\Rab)$. Let $g=\pi^{-1}(\overline{g})\in \varinjlim_{\alpha\in \Laml} \Hom_{\ARing_S}(\hB,R_\lambda)$.
By the commutativity of the diagram. we have
$\iota(p(g))=q(\psi)=\phi=\iota(\tilde{f})$, where $\tilde{f}$ is the image of $f$ in 
$\varinjlim_{\alpha\in \Laml} \Hom_{\ARing_S}(A,\Ra)$. This implies $\tilde{f}=p(g)$ by the injectivity of $\iota$, Remark \ref{rem2;lem;Phi},
which completes the proof.

\begin{claim}\label{claim;lem;compactness}
The natural map
\[\varinjlim_{\alpha\in \Laml} \Hom_{\pARing_A}(B,\Rab)\to \Hom_{\pARing_A}(B,\Rinf)\]
is a bijection.
\end{claim}

Consider the commutative diagram
\[\xymatrix{\varinjlim_{\alpha\in \Laml} \Hom_{\pARing_A}(B,\Rab)\ar[r]\ar[d]^{\subset}
& \Hom_{\pARing_A}(B,\Rinf)\ar[d]^{\subset} \\
\varinjlim_{\alpha\in \Laml} \Hom_{\Ring_A}(B,\Rab)\ar[r] &\Hom_{\Ring_A}(B,\Rinf)\\}\]
It is cartesian by Lemma \ref{lem;cartesian}.
Since $B$ is of finite presentation over $A$, the lower horizontal map is a bijection, which implies the claim. 

\end{proof}

\subsection{Formal schemes}

In this paper, a formal scheme means a topologically ringed space $(\fX,\cO_\fX)$ satisfying the following conditions:
\begin{enumerate}
\item[(a)]
Any point of $\fX$ admits an open neighborhood $\fU$ such that the topologically ringed space induced by $\fX$ on $\fU$ is the formal spectrum $\Spf(A)$ \cite[10.1.2]{EGAI} for some $A\in \aRing$. Such $\fU$ is called an affine formal open of $\fX$.
\item[(b)]
$\fX$ admits an ideal of definition of finite type over $\cO_\fX$.
\end{enumerate}
Here, an ideal $\cI$ of definition of $\fX$ is an ideal $\cI\subset \cO_{\fX}$ such that 
any point of $\fX$ admits an affine formal open neighborhood $\Spf(A)$ such that
$\cI_{|\Spf(A)}$ corresponds to an ideal of definition in $A$ (\cite[10.3.1]{EGAI}). 
It gives rise to a scheme
\eq{cXI}{
\fX_\cI=\uSpec(\cO_\fX/\cI)=(\fX,\cO_\fX/\cI)}
whose underlying topological space $|\fX_\cI|$ is equal to that $|\fX|$ of $\fX$.

A morphism of formal schemes $f:\fY\to \fX$ is a morphism of 
topologically ringed spaces such that for any point $y\in \fY$, it induces a local homomorphism $\cO_{\fX,f(y)}\to \cO_{\fY,y}$ (\cite[10.4.5 and 10.1.5]{EGAI}).
Let $\fSch$ denote the category of formal schemes.
 Note that $\fSch$ admits fiber products such that (\cite[10.4.7]{EGAI}).
\eq{fiberproductfSch}{\Spf(B)\times_{\Spf(A)}\Spf(C)=\Spf(B\hotimes_A C).}

We say that $\fX\in \fSch$ is locally Noetherian if 
any point of $\fX$ admits an affine formal open neighborhood $\Spf(A)$ with $A$ Noetherian. If $\fX$ is locally Noetherian, every ideal of definition of $\fX$ is of finite type (\cite[2.1.28]{Ab10}).

\begin{lemm}\label{lem1;fSch}
Let $\fX\in \fSch$ and $\cI\subset \cO_\fX$ be an ideal of definition of finite type.
For integers $n>0$, put 
$\fX_n=\uSpec(\cO_\fX/\cI^n)=(\fX,\cO_\fX/\cI^n)\in \Sch$
Then, 
\[\fX=\varinjlim_{n>0} \fX_n,\]
where the transition map $\fX_m\to \fX_n$ with $m>n$ is induced by the quotient map $\cO_\fX/\cI^m\to \cO_\fX/\cI^n$ and the colimit is taken in the category $\fSch$.
i.e. it has the same topological space as that of $\fX_n$ and is equipped with the sheaf of topological rings which is the limit of the sheaves of discrete rings $\cO_{\fX_n}$. 
\end{lemm}
\begin{proof}
 See \cite[2.1.23]{Ab10}.
\end{proof}

\begin{defi}\label{def;formalCompletion}
Let $\fX\in \fSch$. A \emph{closed subscheme of $\fX$} is a scheme $Z=\Spec(\cO_\fX/\cI_Z)$ for a coherent ideal $\cI_Z\subset\cO_\fX$ which contains an ideal of definition of $\fX$. The \emph{formal completion $\fc \fX Z$ of $\fX$ along $Z$} is 
associated to 
\[ \varinjlim_n \uSpec(\cO_\fX/\cI_Z^n)\]
via Lemma \ref{lem1;fSch}. Note that there is a natural morphism $\fc \fX Z \to \fX$.
\end{defi}

\begin{rmk}\label{rmk;def;formalCompletion}
If $f:\fX'\to \fX$ is a morphism in $\fSch$, then we have (cf. \eqref{cXI})
\eq{eq;formalCompletion}{
\fc \fX Z\times_{\fX} \fX'=\fc {\fX'}{Z'}\qwith Z'=Z\times_{\fX_\cI}\fX'_{\cI'}, }
where $\cI$ (resp. $\cI'$) is an ideal of definition of $\fX$ (resp. $\fX'$) such that $\cI\subset \cI_Z$ and $f^*(\cI)\subset \cI'$.
This is checked by noting that the fiber products in $\fSch$ are locally given by the formula  \eqref{fiberproductfSch}.
\end{rmk}

\begin{defi}\label{def;adicmorphism}
A morphism $f:\fY\to \fX$ in $\fSch$ is adic (or $\fY$ is adic over $\fX$) if for any ideal $\cI$ of definition of $\fX$, $f^*(\cI)\cO_\fY$ is an ideal of definition of $\fY$. 
Let $\fSchad\subset \fSch$ denote the subcategory of formal schemes with adic morphisms. 
Note that $\fSchad$ is closed under fiber products in $\fSch$.

Let $\cI$ be an ideal of definition of $\fX$ and $\cJ=f^*(\cI)\cO_\fY$.
Then, $f$ induces morphisms in $\SchS$:
\eq{fI}{f_\cI:\fY_\cJ\to \fX_\cI.}
\end{defi}

\begin{rema}
A typical morphism which is not adic is the natural map $\fc \fX Z\to \fX$, where 
$\fc \fX Z$ is the completion of $\fX$ along a closed subscheme 
$Z\subset \uSpec(\cO_\fX/\cI)$ for an ideal $\cI$ of definition of $\fX$, see Definition \ref{def;formalCompletion}.
\end{rema}

\begin{lem}\label{lem1;adicmorphism}
Let $f:\fY\to \fX$ and $g:\fY\to \fZ$ be morphisms in $\fSch$.
\begin{itemize}
\item[(1)]
If $f$ and $g$ are adic, then so is $g\circ f$. 
\item[(2)]
If $g$ and $g\circ f$ are adic, then so is $f$.
\end{itemize}
\end{lem}
\begin{proof}
See \cite[2.2.11]{Ab10}.
\end{proof}

\begin{defi}\label{def;adicmorphism2}
For $\fS\in \fSch$, let $\fSchadfS$ denote the category of adic formal schemes over $\fS$. By Lemma \ref{lem1;adicmorphism}(2), any morphism in $\fSchadfS$ is adic. 
Note that $\fSchadfS$ is closed under fiber products in $\fSchfS$.
\end{defi}

\begin{lem}\label{lem2;adicmorphism}
Let $\fS\in \fSch$ and $\cI$ be an ideal of definition of finite type of $\fS$.
For integers $n>0$, put $\fS_n=(\fS, \cO_{\fS}/\cI^n)$. Then, there is a natural equivalence between $\fSchadfS$ and the category of inductive systems $\{X_n\}_{n>0}$ of $\fS_n$-schemes such that for $m\leq n$, the diagram
\[\xymatrix{
X_m \ar[r]\ar[d] & \fS_m \ar[d]\\
X_n \ar[r]& \fS_n \\}\]
is cartesian, and morphisms are systems of morphisms $u_n:X_n\to Y_n$ of $\fS_n$-schemes such that $u_m=u_n\times_{\fS_n} \id_{\fS_m}$ for $m\leq n$.
\end{lem}
\begin{proof}
See \cite[2.2.14]{Ab10}.
\end{proof}

\subsection{Morphisms of finite type and formally of finite type}

\begin{defi}\label{def;ftfSch}
Let $\fX$ be a locally Noetherian formal scheme and $f:\fY\to \fX$ be a morphism in $\fSch$.
We say that $f$ is \emph{locally of finite type} (resp. \emph{locally formally of finite type}) if for any point $y\in \fY$, there exists an affine formal open neighborhood $U=\Spf(A)$ of $f(y)$ in $\fX$
and an affine formal open neighborhood $V=\Spf(B)$ of $y$ in $U\times_\fX \fY$ such that $B$ is topologically of finite type over $A$ (resp. formally of finite type over $A$), Definition \ref{def;ftAring}.
We say that $f$ is \emph{of finite type} (resp. \emph{formally of finite type}) if it is   quasi-compact and locally of finite type (resp. locally formally of finite type).

\end{defi}

\begin{lemm}\label{lem1;ftfSch}
Let $f:\fY\to \fX$ be as in Definition \ref{def;ftfSch}.
\begin{itemize}
\item[(1)]
$f$ is locally of finite type (resp. of finite type) if and only if $f$ is adic and for some ideal $\cI$ of definition of $\fX$, the morphism $\uSpec(\cO_\fY/\cI\cO_\fY) \to  \uSpec(\cO_\fX/\cI)$ is locally of finite type (resp. of finite type) as a morphism of schemes; then this property holds for any ideal of definition of $\fX$. 
\item[(2)]
$f$ is locally formally of finite type (resp. formally of finite type) if and only if there exists an ideal $\cI$ (resp. $\cJ$) of definition of $\fX$ (resp. $\fY$) such that $\cI\cO_\fY\subset \cJ$ and the morphism $\uSpec(\cO_\fY/\cJ) \to  \uSpec(\cO_\fX/\cI)$ is locally of finite type (resp. of finite type) as a morphism of schemes; then this property holds for any an ideal $\cI$ (resp. $\cJ$) of definition of $\fX$ (resp. $\fY$) such that $\cI\cO_\fY\subset \cJ$. 
\item[(3)]
If $f$ is locally formally of finite type (resp. formally of finite type), it is 
 locally of finite type (resp. of finite type) if and only if it is adic.
\end{itemize}
\end{lemm}
\begin{proof}
This follows from \cite[Pr.1.7]{TLP}, see also \cite[Appendix 3.2.4 and 3.2.5]{CLNS}.
\end{proof}

\begin{coro}\label{cor.lem1;ftfSch}
Morphisms of finite type (resp. formally of finite type) are stable under compositions and base changes, and any morphism between formal schemes of finite type (resp. formally of finite type) over a fixed formal scheme is of finite type (resp. formally of finite type).
\end{coro}
\begin{proof}
This follows from Lemma \ref{lem1;ftfSch} and the corresponding facts for morphisms of finite type for schemes.
\end{proof}


\begin{exam}
\begin{itemize}
\item[(1)]
Let $A$ be a Noetherian adic ring and $I\subset A$ be an ideal of definition and $X$ be a scheme locally of finite type (resp. of finite type) over $A$.
Then, the formal completion of $X$ along $X\times_{\Spec(A)} \Spec(A/I)$ is a formal scheme locally of finite type (resp. of finite type) over $\Spf(A)$.
\item[(2)]
Let $\fX$ be a locally Noetherian formal scheme and $\fc \fX Z$ be the formal completion of $\fX$ along a closed subscheme $Z$, Definition \ref{def;formalCompletion}.
Then, the natural map $\fc \fX Z \to \fX$ is formally of finite type. 
\end{itemize}
\end{exam}

\begin{defi}\label{def;fSchfpS}
For a locally notherian formal scheme $\fS$, let $\fSchfS$ be the category of formal schemes over $\fS$ and 
$\fSchfpfS$ (resp. $\fSchffpfS$) be the full subcategories of $\fSchfS$ of formal schemes of finite type (resp. formally of finite type) over $\fS$. 
\end{defi}

\begin{rmk}\label{rmk;def;fSchfpS}
By Corollary \ref{cor.lem1;ftfSch}, every morphism in $\fSchfpfS$ (resp. $\fSchffpfS$) is of finite type (resp. formally of finite type) and $\fSchfpfS$ (resp. $\fSchffpfS$) is stable under fiber products in $\fSchfS$. 
\end{rmk}

\subsection{\'Etale morphisms and Nisnevich sites}

In this section, all formal schemes are assumed to be locally Noetherian.

\begin{defi}\label{def;etalemorphism and Nisnevich sites}
A morphism $f:\fY\to \fX$ in $\fSch$ is formally \'etale if for any closed immersion $Y_0\to Y$ of schemes defined by a nilpotent ideal and any morphism $Y\to \fX$ in $\fSch$, the induced map
\[ \Hom_{\fSch_\fX}(Y,\fY) \to \Hom_{\fSch_\fX}(Y_0,\fY)\]
is bijective.
The morphism $f$ is \'etale if it is adic, locally of finite type
and formally \'etale. An adic morphism $A\to B$ in $\aRing$ is called \'etale if $\Spf(B)\to \Spf(A)$ is \'etale.
\end{defi}

\begin{lemm}\label{lem;etalemorphism}
\begin{itemize}
\item[(1)]
Let $f:\fY\to \fX$ be an adic morphism in $\fSch$. 
$f$ is \'etale, if and only if the morphism $f_\cI$ from \eqref{fI} is \'etale for every ideal $\cI\subset \cO_\fX$ of definition of finite type.
\item[(2)]
Let $\fX\in \fSch$ with an deal $\cI\subset \cO_\fX$ of definition of finite type.
For an \'etale morphism  $f_0: Y_0\to X_0=\uSpec(\cO_\fX/\cI)$, there exists an \'etale morphism $f:\fY\to \fX$ in $\fSch$ which lifts $f_0$.
Moreover, $f$ is unique up to a unique isomorphism.
\end{itemize}
\end{lemm}
\begin{proof}
See \cite[2.4.8 and 2.4.11]{Ab10}.
\end{proof}

The following deep result on algebrization of \'etale morphisms of adic rings plays an important role in the proof of Theorem \ref{prop:procdhLocal}.

\begin{theo}\label{thm;algebrizaitonEtale}
Let $A$ be a Noetherian adic ring and $I$ be an ideal of definition of $A$.
Let $A\to B$ be an \'etale morphism in $\ARing$.
Then, there exists an \'etale $A$-algebra $B'$ whose $IB'$-adic completion is isomorphic to $B$. 
\end{theo}
\begin{proof}
This is a consequence of a more general result \cite[III th\'eo. 7]{El}, 
see \cite[1.16.27 and 1.16.29]{Ab10}.
\end{proof}

\begin{defi}\label{def;fXnis}
For $\fX\in \fSch$, let $\fX_{\Nis}$ be the following site: 
The underlying category is the full subcategory of $\fSch_{\fX}$ of \'etale morphisms $\fY\to \fX$.
A sieve on $\fY$ is a covering if it contains a finite family
$ \{\fY_i \to \fY\}_{i\in I}$ such that 
\begin{enumerate}
\item[$(\spadesuit)$]
There exists an ideal $\cJ$ of definition of finite type of $\fY$ such that 
\[ \{(\fY_i,\cO_{\fY_i}/\cJ\cO_{\fY_i}) \to (\fY,\cO_{\fY}/\cJ)\}_{i\in I}\]
is a Nisnevich covering.
\end{enumerate}
\end{defi}

\begin{rmk}\label{rmk;def;fXnis}
If the condition $(\spadesuit)$ holds for some $\cJ$, it holds any ideal $\cJ$ of definition of finite type of $\fY$. Hence, by Lemma \ref{lem;etalemorphism},
for any ideal $\cI$ of definition of finite type of $\fX$, the restriction induces an equivalence of sites
\eq{Nisequiv}{
\fX_\Nis \simeq (\fX_\cI)_\Nis\qwith \fX_\cI=\uSpec(\cO_\fX/\cI))=(\fX,\cO_\fX/\cI),}
where the right hand side is the usual Nisnevich site on the scheme $\fX_\cI$.
\end{rmk}

\begin{defi} \label{defi:distinguishedNissquare}
A \emph{distinguished Nisnevich square} is a square in $\fX_\Nis$:
\eq{Nissquare}{\xymatrix{
\fY' \times_\fY \fU \ar[r] \ar[d] & \fY'\ar[d]^f \\
\fU \ar[r]^j & \fY\\}}
where $j$ is an open immersion of formal schemes and that $f$ is an \'etale morphism of affine formal schemes over $\fX$ such that the induced morphism 
$\fY'\times_\fX \fX_\cI \to \fY\times_\fX \fX_\cI$ of schemes for an ideal $\cI$ of definition of $\fX$ is an isomorphism over a closed subscheme $Z$ of $\fY\times_\fX \fX_\cI$ such that 
$|Z|=|\fY|\setminus |\fU|$.
\end{defi}

\begin{prop}\label{prop;NisfX}
The Nisnevich topology on $\fX$ is generated by the Zariski topology and the covering family $\fU\sqcup \fY' \to \fY$ for distinguished Nisnevich squares \eqref{Nissquare}.
\end{prop}
\begin{proof}
In view of \eqref{Nisequiv}, this follows from the corresponding fact on the Nisnevich topology on $\fX_\cI$
\end{proof}

\subsection{Blowups}\label{Blowup}

\begin{defi}\label{defi;Blowup}
Let $\fX\in \fSch$ and take an ideal $\cI$ of definition of finite type of $\fX$ and
$\fX_n=(\fX, \cO_{\fX}/\cI^n)$ for integers $n>0$.
For a coherent open ideal $\cA\subset \cO_\fX$, we put
\[ \fX_n'=\Proj(\underset{r\geq 0}{\bigoplus} \cA^r\otimes_{\cO_\fX} \cO_{\fX_n}).\]
For $m\leq n$, the canonical morphism $\fX_m\to \fX_n$ induces a morphism
$\fX'_m\to \fX'_n$ such that the diagram
\[\xymatrix{
\fX'_m \ar[r]\ar[d] & \fX_m \ar[d]\\
\fX'_n \ar[r]& \fX_n \\}\]
is cartesian, where the vertical maps are the natural projections.
By Lemma \ref{lem2;adicmorphism}, we get an adic formal scheme $\fX'$ over $\fX$ as the inductive limit of $\fX'_n$. It follows directly from the definition that the natural projection $\fX'\to \fX$ is proper and of finite type.
We call it the blowup of $\cA$ in $\fX$, or the blowup of $\fX$ in the closed subscheme $\uSpec(\cO_{\fX}/\cA)\subset\fX$.
\end{defi}

\begin{lem}\label{lem1;Blowup}
Let $X$ be a scheme, $Y$ be a closed subscheme of $X$ defined by a coherent ideal of $\cO_X$ and $U=X\setminus Y$ and $\fX=\fc X Y$ be the formal completion of $X$ along $Y$. 
Let $\cA\subset \cO_X$ be a coherent ideal such that $\cA_{|U}=\cO_U$ and $\phi: X'\to X$ be the blowup of $\cA$ in $X$.
Then, the blowup of $\cA_{|Y}$ in $\fX$, Definition \ref{defi;Blowup},
 is canonically isomorphic to the formal completion $\fc {X'} {Y'}$ of $X'$ along $Y'=Y\times_X X'$.
\end{lem}
\begin{proof}
This is \cite[3.1.3]{Ab10}.
\end{proof} 

\begin{defi}\label{strict-peoprTransform}
Let $f:\fX'\to \fX$ be the blowup of $\fX\in \fSch$ in a closed subscheme 
$\Sigma\subset\fX$, Definitions \ref{def;formalCompletion} and \ref{defi;Blowup}.
Let $T$ be a closed subscheme of $\fX$ and take an ideal $\cI$ of definition of $\fX$ such that $\Sigma$ and $T$ are closed subschemes of $X_0=\uSpec(\cO_\fX/\cI)$.
Then, $f$ is an isomorphism over $U=X_0\setminus \Sigma$ so that $U$ is viewed as an open subscheme of $X_0'=X_0\times_\fX\fX'$. Then, the \emph{proper transform of $T$ in $\fX'$} is the closed subscheme of $\fX'$ which is defined as the closure of $T\cap U$ in $X_0'$.
\end{defi}



\section{Definition of a pro-cdh topology}

In this section, we fix a locally Noetherian formal scheme $S$
and define a pro-cdh topology on $\fSchffpS$.

\begin{defi} \label{defi:locallyalgebriaable}
Let $\fX$ be a locally Noetherian formal scheme and $Z=\uSpec(\cO_\fX/\cI_Z)$ is a closed subscheme of $\fX$, Definition \ref{def;formalCompletion}. 
We let $\LAP(\fX,Z)$ denote the class of adic morphisms $f:\fY\to \fX$ satisfying the conditions:
\begin{enumerate}
\item[(i)]
$f$ is proper, 
i.e. the map $f_\cI:\fY_\cJ\to \fX_\cI$ of schemes from \eqref{fI} is proper for an ideal $\cI$ of definition of $\fX$, 
and of finite type, Definition \ref{def;ftfSch}.
\item[(ii)]
$f$ is \emph{$Z^c$-admissible}, i.e. an isomorphism over the open complement of $Z$ in $\fX$.
\item[(iii)]
$f$ is \emph{locally algebrizable}, i.e. any point of $\fX$ admits an affine formal open neighborhood $\Spf(A)$ such that 
there is a morphism $g:Y\to \Spec(A)$ of schemes which is proper of finite type such that 
\[\fY\times_\fX \Spf(A)=Y\times_{\Spec(A)}\Spf(A).\]
Moreover, $g$ is an isomorphism over $\Spec(A)\setminus T$ for some closed subscheme $T$ such that $T\cap \Spec(A/I)\subset Z$ for an ideal $I\subset A$ of definition.

\end{enumerate}
\end{defi}

\begin{rema} \label{rema1:locallyalgebrisable}
By Lemma \ref{lem1;Blowup} (also \cite[3.1.1.3]{Ab10}), if $f:\fY\to \fX$ is the blowup of a coherent open ideal $\cA\subset \cO_\fX$ in $\fX$, then $f$ is in $\LAP(\fX,Z)$ with $Z=\uSpec(\cO_X/\cA)$. 
\end{rema}

\begin{lemm}\label{lem0:locallyalgebrisable}
Let $f:\fY\to \fX$ be in $\LAP(\fX,Z)$.
\begin{itemize}
\item[(1)]
For any closed subscheme $Z'$ of $\fX$ such that $|Z|\subset |Z'|$, $f$ is in $\LAP(\fX,Z')$.
\item[(2)]
For $g:\fX'\to \fX$ in $\fSch$, $\fY\times_\fX \fX'\to \fX'$ is in $\LAP(\fX',Z')$, where choosing an ideal $\cI$ (resp. $\cI'$) of definition of $\fX$ (resp. $\fX'$) such that $\cI\subset \cI_Z$ and $g^*(\cI)\subset \cI'$, $Z'=Z\times_{\fX_\cI} \fX'_{\cI'}$ (cf, \eqref{cXI}).
\item[(3)]
The diagonal $\Delta_\fY:\fY \to \fY\times_\fX \fY$ is in $\LAP(\fY\times_\fX \fY,W\times_Z W)$, where choosing an ideal $\cI$ of definition of $\fX$ such that $\cI\subset \cI_Z$, $W=Z\times_{\fX_\cI} \fY_{\cJ}$ with $\cJ=\cI\cO_\fY$.
\end{itemize}
\end{lemm}
\begin{proof}
(1) is obvious. We prove (2). Let the notation be as in Definition \ref{defi:locallyalgebriaable}(iii). We may assume $\fX=\Spf(A)$ and $\fX'=\Spf(A')$ with ideals of definition $I\subset A$ and $I'\subset A'$ such that $IA'\subset I'$ and 
$I\subset J$ for the ideal $J$ defining $Z$.
Putting $g':Y'=Y\times_{\Spec(A)}\Spec(A')\to \Spec(A')$, we have 
$\fY\times_\fX \fX'=Y'\times_{\Spec(A')}\Spf(A')$ and $g'$ is an isomorphism over
$\Spec(A')\setminus T'$, where $T'=T\times_{\Spec(A)}\Spec(A')$.
Letting $K\subset A$ be the ideal defining $T$, 
the assumption $T\cap \Spec(A/I)\subset Z$ implies $J\subset I+K$ so that 
$I'+JA'\subset I'+KA'$, and hence, $T'\cap \Spec(A'/I')\subset Z'$, where 
$Z'=Z\times_{\Spec(A/I)}\Spec(A'/I')$. This proves (2).

To prove (3), we may assume $\fX=\Spf(A)$. Consider the following diagram where the squares are carteisan in $\fSch$:
\[\xymatrix{
\fY \ar[r]^-{\Delta_\fY}\ar[d] & \fY\times_\fX \fY \ar[r]\ar[d] & \Spf(A) \ar[d]\\
Y \ar[r]^-{\Delta_Y} & Y\times_{\Spec(A)} Y \ar[r]^-f & \Spec(A) \\}\]
The maps $\Delta_Y$ and $f$ are isomorphisms over $\Spec(A)\setminus T$ and
the assumption $T\cap \Spec(A/I)\subset Z$ implies
$f^{-1}(T)\cap f^{-1}(\Spec(A/I)) \subset f^{-1}(Z)=W\times_Z W$, which proves (3).


\end{proof}

\begin{defi} \label{defi:procdh}
The \emph{pro-cdh topology} on $\fSchffpS$ is generated by the following coverings.
\begin{enumerate} 
 \item[(Nis)]  Nisnevich coverings, Definition \ref{def;fXnis}.

\item[(FAB)] Formal abstract blowup covering: families of the form
 \[ \{\fc \fX Z \sqcup \fY \to \fX\} \]
 where $Z$ is a closed subscheme of $\fX$ and $\fY \to \fX$ is in $\LAP(\fX,Z)$.
\end{enumerate}

We let $\Shvzen(\fSchffpS)$ denote the category of pro-cdh sheaves of sets.
\end{defi}

\begin{rema}\label{rema3:locallyalgebrisable}
In view of Lemma \ref{lem0:locallyalgebrisable}(2) and \eqref{eq;formalCompletion},
the covering family (FAB) is stable under base changes.
\end{rema}

\begin{rema}\label{rem;defi:procdh}
Let $\SchftSb$ be the full subcategory of the category of schemes of finite type over $S$ consisting those $X\to S$ that factor through $\uSpec(\cO_S/I)$ for some ideal $I$ of definition of $S$. 
For such $X\in\SchftSb$, the map of $X\to \Spec(S/I)\hookrightarrow \Spf(S)$ is an object of $\fSchffpS$ so that $\SchftSb$ is viewed as a full subcategory of $\fSchffpS$.
The embedding $\iota:\SchftSb\to \fSchffpS$ induces a pair of adjoint functors
\begin{equation}\label{adjunction-iota}
\cont: \PSh(\SchftSb)\begin{smallmatrix}\longleftarrow \\ \longrightarrow\\
\end{smallmatrix} \PSh(\fSchffpS):\iota^* 
\end{equation}
where $\iota^*$ (resp. $\cont$) is the restriction (resp. the right Kan extension) along $\iota$. It induces an equivalence
\eq{PShprocdhcont}{\cont: \PSh(\SchftSb) \isom \PSh^{\cont}(\fSchffpS),}
where $\PSh^{\cont}(\fSchffpS)$ denotes the full subcategory of those presheaves satisfying 
\begin{equation} \label{equa:cont}
F(\fX) = \varprojlim_n F(\fX_n)\qfor \fX\in \fSchffpS,
\end{equation}
where $\fX_n=\uSpec(\cO_\fX/\cI^n)$ for an ideal $\cI$ of definition of $\fX$.
The continuity condition \eqref{equa:cont} is equivalent to the sheaf condition for the topology on $\fSchffpS$ whose covering families are those of the form
\begin{enumerate}
 \item[(Cont)] $\{ \fX_n \to \fX \}_{n \in \NN}$ for $\fX\in\fSchffpS$ .
\end{enumerate}
The equivalence can be seen by noting $\fX_n \times_{\fX}\fX_m = \fX_{\min(n,m)}$ and 
\[\varprojlim_n F(\fX_n) = eq (\prod_\NN  F(\fX_n) \rightrightarrows \prod_\NN  F(\fX_n)).\]
One sees that \eqref{PShprocdhcont} restricts to an equivalence
\eq{Shvprocdhcont}{
\cont: \Shvzen(\SchftSb) \isom \Shvzenc(\fSchffpS),}
where $\Shvzen(\SchftSb)$ is the full subcategory of $\PSh(\SchftSb)$ 
consisting of those presheaves whose restrictions to each $\Schft_{S/I^n}$ for $n>0$
is a pro-cdh sheaf in the sense of \cite[Def. 1.1]{KS23}
and $\Shvzenc(\fSchffpS)$ is the category of sheaves on $\fSchffpS$ for the topology generated by (Nis) and (FAB) and (Cont). 
\end{rema}

In what follows, we address some basic properties of the pro-cdh topology.

\begin{defi}\label{def:closedcovering}
Let $\fX\in \fSchffpS$ and $Z,W$ be closed subschemes of $\fX$, Definition \ref{def;formalCompletion}, such that $|\fX|=|Z|\cup |W|$.
Then, the square
\eq{ClosedCoveringSquare} 
{ \xymatrix{ 
\fc {\fX} {Z\cap W} \ar[r] \ar[d] & \fc {\fX} {W} \ar[d]\\
\fc {\fX} {Z} \ar[r]  & \fX \\}  }
is called a closed coverings square in $\fSchffpS$, where $Z\cap W=Z\times_\fX W$.
\end{defi}

\begin{lemm}\label{lem;ClosedCoveringSquaire}
For every square \eqref{ClosedCoveringSquare}, $\fc {\fX} {Z} \sqcup \fc {\fX} {W} \to\fX$ is a pro-cdh covering.  
\end{lemm}
\begin{proof}
The proof is suggested by J. Ayoub.
Let $\fX'=\Bl_{Z\cap W}(\fX)$ be the blowup of $\fX$ in $Z\cap W$, Definition \ref{defi;Blowup} and $Z'$ (resp. $W'$) be the proper transform of $Z$ (resp. $W$) in $\fX'$, Definition \ref{strict-peoprTransform} and put 
$\fU'=\fX\setminus Z'$ and $\fV'=\fX\setminus W'$.
Then, we have $Z'\cap W'=\emptyset$ so that $\fU'\sqcup \fV'\to \fX'$ is a Zariski covering. Composing this with a pro-cdh covering $\fc {\fX} {Z\cap W}\sqcup \fX' \to \fX$, we get a pro-cdh covering $\cU'\sqcup\cV'\sqcup \fc{\fX} {Z\cup W}\to \fX$.
Since $\fU' \to \fX$ (resp. $\fV'\to \fX$) factors through $\fc\fX W$ (resp. $\fc\fX Z$), this proves the lemma.

\end{proof}

The following lemma was suggested by J. Ayoub.

\begin{lemm}\label{lem;Joseph}
The pro-cdh topology on $\fSchffpS$ is generated by Nisnevich topology and by families of the form
\eq{eq;lem;Joseph}{  \fc \fX Z\sqcup \fc {(\Bl_\Sigma \fX)} {X'_0} \to \fX,}
where $\fX\in \fSchffpS$ and 
$\Sigma\subset Z\subset X_0:=\uSpec(\cO_\fX/\cI)$ are closed subschemes 
for some ideal $\cI\subset\cO_\fX$ of definition and  $\Bl_\Sigma \fX\to \fX$ is the blowup of $\fX$ in $\Sigma$, 
Definition \ref{defi;Blowup}, and $\fc {(\Bl_\Sigma \fX)} {X'_0}$ is the formal completion of 
$\Bl_\Sigma \fX$ along the proper transform $X_0'$ of $X_0$, Definition \ref{strict-peoprTransform}.
\end{lemm}
\begin{proof}(Ayoub)
First, we claim that \eqref{eq;lem;Joseph} is a covering. Indeed, 
$\fc \fX Z\sqcup \fY\to \fX$ with $\fY=\Bl_\Sigma \fX$ is a covering by Definition \ref{defi:procdh} and $\fc\fY {X_0'}\sqcup \fc \fY E \to \fY$ with $E=\fY\times_{\fX} \Sigma$ is a covering by Lemma \ref{lem;ClosedCoveringSquaire}. Hence, the claim follows from the fact that 
$\fc \fY E \to \fX$ factors through $\fc \fX \Sigma\to \fc \fX Z$.

It remains to prove that for $\fY \to \fX$ in $\LAP(\fX,Z)$ as in Definition \ref{defi:procdh}(FAB), we can refine the covering $\fc \fX Z\sqcup \fY \to \fX$ by a covering of the form
\eqref{lem;Joseph}. By Definition \ref{defi:locallyalgebriaable}, we may assume that $\fX=\Spf(A)$ and there is a morphism $g:Y\to X=\Spec(A)$ of schemes which is proper of finite type such that 
\[\fY\times_\fX \Spf(A)=Y\times_{\Spec(A)}\Spf(A),\]
and that $g$ is an isomorphism over $\Spec(A)\setminus T$ for some closed subscheme $T$ such that $T\cap \Spec(A/I)\subset Z$ with $I=\cI A$.
By \cite[Se.5]{RG71} and \cite[Lem.2.1.5]{Tem08}, there exist a closed subshcme $T'\subset T$ such that 
$\Bl_{T'}(X)\to X$ factors through $Y$. Hence, we may assume $Y=\Bl_T(X)$.
Then, letting $f:Y'=\Bl_\Sigma(X)\to X$ with $\Sigma:=T\cap \Spec(A/I)$ and $X_0'$ be the proper transform of $X_0$ in $Y'$,
it suffices to show that $\fc{Y'} {X_0'}\to X$ factors through $\Bl_T(X)$.
Let $T'$ be the proper transform of $T$ in $Y'$ and $E=f^{-1}(\Sigma)$.
Since the ideal sheaf $\cI_{T'}$ is equal to $\cI_{f^{-1}(T)}\cdot (\cI_E)^{-1}$,
we have a map $\Bl_{T'}(Y') \to \Bl_T(Y)$ by the universality of blowups. 
On the other hand, we have a map $\fc{Y'} {X_0'} \to \Bl_{T'}(Y') $ since $X_0'\cap T'=\emptyset$. This proves the desired assertion.
 \end{proof}

\begin{lemm} \label{lemm:squareDescent}
Let $F\in \PSh(\fSchffpS,\Spc)$ be a presheaf of spaces on $\fSchffpS$.
\begin{itemize}
\item[(1)]
$F$ satisfies {\v Cech} descent for the family of coverings of the form $\fY\sqcup \fc\fX Z \to \fX$ from Definition \ref{defi:procdh} if and only if it sends all squares 
\eq{cdh-square-h}{\xymatrix{
\fc \fX Z\times_\fX\fY \ar[r]^-{j}\ar[d] & \fY \ar[d]^{f} \\
\fc \fX Z \ar[r]^-{i} & \fX\\}}
to cartesian squares. 
\item[(2)]
$F$ satisfies {\v Cech} descent for the family of coverings of the form $\fc {\fX} {Z} \sqcup \fc {\fX} {W} \to\fX$ from Lemma \ref{lem;ClosedCoveringSquaire}) if and only if it sends all squares \eqref{ClosedCoveringSquare} to cartesian squares. 
\end{itemize}
\end{lemm}
\begin{proof}
Consider the cd-structure on $\fSchffpS$ associated to all squares of the form \eqref{cdh-square-h}.
By Lemma \ref{lem0:locallyalgebrisable}, those squares are stable under base changes and the squares
\[\xymatrix{ 
\fc {\fX} {Z}\times_\fX \fY \ar[r] \ar[d] & \fY \ar[d]\\
\fc {\fX} {Z}\times_\fX\fY\times_\fX \fY\ar[r]  
& \fY\times_\fX \fY\\} \]
with the diagonal vertical maps are again of the form \eqref{cdh-square-h}.
Hence, (1) follows from \cite[Lem. 2.5 and 2.11]{VoeHTS} and \cite[Th.3.2.5]{AHW17}.
(2) is proved by the same argument noting that for a given square \eqref{ClosedCoveringSquare}, the square
\[\xymatrix{ 
\fc {\fX} {Z\cap W} \ar[r] \ar[d] & \fc \fX W\ar[d]\\
\fc {\fX} {Z\cap W}\times_{\fc {\fX} {Z}} \fc {\fX} {Z\cap W} \ar[r]  
& \fc {\fX} {W}\times_\fX \fc {\fX} {W}\\} \]
with the diagonal vertical maps are again of the form \eqref{ClosedCoveringSquare} 
noting that the vertical maps are isomorphisms.

\end{proof}


\section{Homotopy dimension}

The main result of this section is the following.

\begin{theo}\label{thm2;Hypercomplete}
Let $S$ be a locally Noetherian formal scheme. Take $\fX \in \fSchffpS$ and let $d$ be the Krull dimension of $|\fX|$.
Then, the $\infty$-topos $\Shvzen(\fSchffpSX,\cS)$ has homotopy dimension$\leq d$ (\cite[Prop.6.5.1.12, Def.7.2.1.1]{HTT}).
\end{theo}


The proof of the above theorem follows the strategy in \cite[\S7]{KS23}.
For this, we use Nisnevich-Riemann-Zariski spaces of  $\fX\in \fSchffpS$, which are defined as 
\[ \RZ(\fX_\Nis) = \varinjlim_{\fY\in \Mod(\fX)} \fY_{\Nis},\]
with the colimit being indexed by modifications $\fY \to \fX$ and taking place in the category of finitary ∞-sites (see Definition \ref{def;RZ}).


\subsection{Modifications} \label{sec:Modifications}


\begin{defi}\label{def;Modificcation}
Let $\fX$ be a locally Noetherian formal scheme and $U\subset \fX$ be a quasi-compact open.
A morphism in $\fSchffpX$ is \emph{$U$-admissible} if it is a composite of morphisms of the following types:
\begin{enumerate}
\item[(i)]
a morphism $f:\fY' \to \fY$ is in $\LAP(\fY,T)$, Definition \ref{defi:locallyalgebriaable}, for some closed subscheme $T$ of $\fY$, Definition \ref{def;formalCompletion}, such that $|\fY\setminus T| \supset |U\times_\fX\fY|$.
\item[(ii)]
a morphism $\fc \fY W\to \fY$, where $W$ is a closed subscheme of $\fY$ such that $W \supset |U\times_\fX\fY|$.
\end{enumerate}
An object $\fY\in \fSchffpX$ is called a $U$-modification if $\fY\to \fX$ is $U$-admissible. Note that every $U$-modification $\fY\to \fX$ is an isomorphism over $U$.

The subcategory of $U$-modifications with $U$-admissible morphisms is denoted $\Mod_U(\fX) ⊆ \fSchffp_\fX$. We also define
\[ \Mod(\fX) = \varinjlim_{U\subset \fX} \Mod_U(\fX) ,\]
where $U$ ranges over all open in $\fX$ such that $|U|$ is dense in $|\fX|$.
Note that the colimit is filtered since the intersection of two dense open subsets in $|\fX|$ is dense.

\end{defi}

\begin{rema} \label{rema:ModLimits} 
By Lemma \ref{lem0:locallyalgebrisable} and \eqref{eq;formalCompletion},
if $\fY' \to \fY$ and $\fY'' \to \fY$ are morphisms in $\ModX$, then $\fY' \times_{\fY} \fY''\to \fY$ and the diagonal $\fY'\to \fY'\times_\fY\fY'$ are again in $\ModX$. 
In particular, $\ModX$ admits finite limits, calculated in $\fSchffpX$, and is therefore is filtered.
\end{rema}



\begin{rema} \label{rema:genXY}
For every $f:\fY \to \fX$ in $\Mod_U(\fX)$, there is $\fY' \to \fY$ in $\Mod_U(\fX)$ such that $|U\times_\fX \fY'|$ is dense in $|\fY'|$.
Indeed, take an ideal $\cI$ (resp. $\cJ$) of definition of $\fX$ (resp. $\fY$) such that $f^*(\cI)\subset \cJ$.
Let $U_\cI=U\times_\fX \fX_\cI$ and $Z$ be the closure of $U_\cI\times_{\fX_\cI} \fY_{\cJ}$ in $\fY_{\cJ}$.
Then, we may take $\fY'=\fc \fY Z$. 
\end{rema}

%

\begin{lemm}\label{excisionAdmMod}
Let $f:\fY' \to \fY$ be a morphism in $\Mod_U(\fX)$ for a quasi-compact open $U$ in $\fX$.
Let $Z$ be a closed subscheme of $\fY$ such that $|U\times_\fX\fY| \supset |\fY\setminus Z| $ and put $Z'=Z\times_\fY \fY'$.
Then, every pro-cdh sheaf $F\in \Shvzen(\fSchffpX,\Spc)$ of spaces sends the square
\eq{eq;excisionAdmMod}{
\xymatrix{\fc{\fY'}{Z'} \ar[d]\ar[r] &\fY'\ar[d]^f\\ 
\fc{\fY} Z \ar[r] &\fY\\ }}
to a cartesian square of spaces.
\end{lemm}
\begin{proof}
Note $\fc{\fY'}{Z'}=\fc{\fY} Z\times_\fY\fY'$ by \eqref{eq;formalCompletion} and that $|U\times_\fX\fY| \supset |\fY\setminus Z| $ implies 
$|U\times_\fX\fY'| \supset |\fY'\setminus Z'|$.
Since the compositions of cartesian squares is cartesian, we may assume that $f$ is either of the type (i) or (ii) from Definition \ref{def;Modificcation}.
For the type (i), we have $|\fY\setminus Z|\subset |U\times_\fX\fY| \subset |\fY\setminus T| $ so that $|T|\subset |Z|$. By Lemma \ref{lem0:locallyalgebrisable}(1),
$\fc\fY Z\sqcup \fY' \to \fY$ is a covering (FAB) from Definition \ref{defi:procdh}. 
So, the assertion follows from Lemma \ref{lemm:squareDescent}(1).
For the type (ii), we may write $\fY'=\fc\fY W$ for a closed subscheme $W$ of $\fY$ such that $W \supset |U\times_\fX\fY|$.
Then, we have $|\fY|=|Z|\cup |W|$ and the assertion follows from Lemma \ref{lemm:squareDescent}(2).
\end{proof}

\subsection{Nisnevich-Riemann-Zariski spaces} \label{sec:NRZspaces}

In this subsection we consider a Nisnevich version of the Riemann-Zariski space associated to a locally Noetherian formal scheme $\fX$.
This will be used as small sites. 
The main result of this section is Corollary~\ref{coro:RZColimitLim} which says that for each $\fX \in \fSchffpS$, the canonical comparison functor $\Shv_{\procdh}(\fSchffpS) \to \Shv(\RZ(\fX_{\Nis}))$ preserves colimits and finite limits.
This result plays a crucial role in the proof of Theorem \ref{thm2;Hypercomplete}.

\begin{defi}\label{def;RZ}
For a locally Noetherian formal scheme $\fX$ we define
\[ \RZ(\fX_{\Nis}) = \int_{\fY \in \ModX} \fY_{\Nis}. \]
Explicitly, $\RZ(\fX_{\Nis}) \subseteq \Arr(\fSchffpX)$ is the category whose objects are morphisms $\fU \to \fY$ such that $\fY \in \ModX$ and $\fU \in \fY_{\Nis}$, Definition \ref{def;fXnis}, and morphisms are commutative squares
\[ \xymatrix{
\fU' \ar[d] \ar[r] & \fU \ar[d] \\
\fY' \ar[r] & \fY
} \]
We abbreviate $\fU \to \fY$ to $(\fU/\fY)$.
\end{defi}

\begin{rema}(see \cite[Rem.4.5]{KS23})
$\RZ(\fX_{\Nis})$ admits finite limits, and they are calculated termwise.
\end{rema}

%
%

\begin{defi} \label{defi:RZtop}
We equip $\RZ(\fX_{\Nis})$ with the Grothendieck topology generated by:
\begin{enumerate}
 \item families of the form
\begin{equation} \label{equa:RZNisCov} \tag{Nis}
\{ (\fU_i/\fY) \to (\fU/\fY) \}_{i \in I}	
\end{equation}
such that $\{\fU_i \to \fU\}$ is a Nisnevich covering, Definition \ref{def;fXnis}, 
 
 \item families of the form
\begin{equation} \label{equa:cartCov} \tag{Car}
\{(\fU\times_\fY \fY' /\fY') \to (\fU/\fY)\},
\end{equation}
for $\fY'\to \fY$ in $\ModX$.	
%
\end{enumerate}
\end{defi}

\begin{lemm} \label{lemm:mNisNormalForm}
Any covering in $\RZ(\fX_\Nis)$ is refinable by one of the form
\[ \{(\fU'_i/\fY') \to (\fU'/\fY') \to (\fU/\fY)\}_{i \in I} \]
where $\fY' \to \fY \in \ModX$ and $\fU'=\fU\times_\fY \fY'$ and 
$\{\fU_i \to \fU'\}_{i \in I}$ is a Nisnevich covering.
\end{lemm}
\begin{proof}
It suffices to show that given families $\{(\fU_i/\fY) \to (\fU/\fY)\}_{i\in I}$ of the form \eqref{equa:RZNisCov} and $(\fV_i/\fY_i) \to (\fU_i/\fY)$ ($\fV_i=\fU_i\times_{\fY} \fY_i$) of the form \eqref{equa:cartCov} for each $i\in I$, we can find $\fY'\to \fY$ in $\ModX$ such that 
$\fU_i\times_{\fY} \fY' \to \fY$ factors through $\fV_i \to \fY_i \to \fY$ for each $i\in I$.
Since $\fU$ is quasi-compact, we can assume that $I$ is finite and
let $\fY'$ be the fibre product of all the $\fY_i$ over $\fY$ (cf. Remark \ref{rema:ModLimits}).
Now we form the following commutative diagram
\[\xymatrix{
&&\ar[lld] \fU\times_{\fY}\fY' \\
\fU\ar[d] &\ar[l] \fU_i \ar[d] & \ar[l] \fV_i \ar[d] & 
\ar[lu]\ar[l] \fV_i' \ar[d] \\
\fY &\ar[l] \fY& \ar[l] \fY_i & \ar[l] \fY'\\}\]
where the right and middle squares are cartesian in $\fSchffpS$.
Then, $\{\fV'_i \to \fU\times_{\fY}\fY'\}_{i\in I}$ is a Nisnevich covering. This proves the lemma.
\end{proof}

\begin{rema}\label{rema:carLoc}(cf. \cite[Rem.4.7 and 4.8]{KS23})
Since the diagonal of a modification is again a modification by Remark \ref{rema:ModLimits}, a presheaf of sets (resp. spaces) satisfies descent for all families \eqref{equa:cartCov} if and only if it sends each 
$(\fY' \times_\fY \fU /\fY') \to (\fU/\fY)$ to an isomorphism (resp. equivalence). Consequently, 
\eq{ShvRZ}{ \Shv(\RZ(\fX_{\Nis})) = \varinjlim_{\fY\in \Mod(\fX)} \Shv(\fY_{\Nis}), }
where the limit is along pushforwards $f_*: \Shv(\fY'_{\Nis}) \to \Shv(\fY_{\Nis})$ for morphisms $f: \fY' \to \fY$ in $\ModX$.
\end{rema}

%


The following proposition is proved by the verbatim same argument as the proof of \cite[Pr.4.8]{KS23}.
 
\begin{lemm} \label{prop:sheafificationCalc}
Suppose $F \in \PSh(\RZ(\fX_{\Nis}))$ satisfies descent for the coverings \eqref{equa:RZNisCov}. Then the sheafification $aF \in \Shv(\RZ(\fX_{\Nis}))$ satisfies
\begin{equation} \label{equa:aFUY}
aF(\fU/\fY) = \varinjlim_{\fY' \in (\ModX)_{/\fY}} F(\fY' \times_\fY \fU / \fY').
\end{equation}
The same is true for presheaves of spaces.
\end{lemm}

\begin{defi}\label{def;rho}
Let $S$ be a locally Noetherian formal scheme.
For $\fX\in \fSchffpS$, we consider the canonical projection functor
\[ \rho_\fX: \RZ(\fX_{\Nis}) \to \fSchffpS; \quad (\fU/\fY) \mapsto \fU \]
and the functor induced by composition
\eq{rho}{ \PSh(\fSchffpS) \to \PSh(\RZ(\fX_{\Nis})); \qquad F \mapsto F \circ \rho_\fX. }
By composing this with the sheafification functor $\PSh(\RZ(\fX_{\Nis})) \to \Shv(\RZ(\fX_{\Nis}))$, we get  
\begin{equation}\label{rhoX}
\rho^*_\fX: \Shv_{\procdh}(\fSchffpS) \stackrel{}{\to} \Shv(\RZ(\fX_{\Nis})). 
\end{equation}
\end{defi}

Recall that a morphism of sites $\phi: C \to D$ is \emph{cocontinuous} if for every $X \in C$ and covering family $\sU = \{U_i \to \phi X \}_{i \in I}$ in $D$, there is a covering family $\{V_i \to X\}_{i \in J}$ in $C$ such that $\{\phi V_i \to \phi X\}_{i \in J}$ refines $\sU$, \cite[Def.III.2.1]{SGA41}, \cite[00XJ]{stacks-project}.

\begin{prop} \label{prop:cocont}
$\rho_\fX$ is cocontinuous. 
\end{prop}

\begin{proof}
For $(\fU/\fY) \in \RZ(\fX_{\Nis})$ and a pro-cdh covering $\fV \to \fU$ in $\fSchffpS$, 
we want to find a morphism $\fY'\to \fY$ in $\ModX$ and a Nisnevich covering:
\eq{equa:VZUVWU0}{ \{ \fU_i\to \fU\times_{\fY} \fY'\}_{i\in I}}
and a map $\phi: \coprod_{i\in I} \fU_i \to \fV$
fitting into a commutative diagram
\eq{equa:VZUVWU}{ \xymatrix{ \coprod_{i\in I} \fU_i  \ar[r]\ar[d]^\phi & \fU\times_{\fY} \fY' \ar[d] \\ \fV \ar[r] & \fU\\}}
Since pro-cdh coverings are refined by finite length compositions of generator pro-cdh coverings, 
it suffices to prove the claim in case $\fV\to \fU$ is a Nisnevich covering, Definition \ref{def;fXnis} or a covering from Lemma \ref{lem;Joseph}.

\def\UI{U_0}
\def\tUI{\tU_0}
\def\fUI{U_0}
\def\tfUI{\tilde{U}_0}
\def\fYI{Y_0}

For a Nisnevich covering the statement is obvious since for any $(\fU/\fY) \in \RZ(\fX_{\Nis})$, a Nisnevich covering $\{\fU_i \to \fU\}_{i \in I}$ of $\fU$ gives rise to a Nisnevich covering $\{(\fU_i/\fY) \to (\fU/\fY)\}_{i \in I}$ of $(\fU/\fY)$ in $\RZ(\fX_\Nis)$. 

Next we treat a covering from \eqref{eq;lem;Joseph}
 \[ \fc \fU {\Xi'} \sqcup \fc{\fW} {W_0} \to \fU \qwith \fW=\Bl_\Xi(\fU),\]
where $\Xi\subset \Xi'$ are closed subschemes of $U_0=Y_0\times_\fY \fU$
with $Y_0=\uSpec(\cO_\fY/\cI)$ for an ideal $\cI$ of definition of $\fY$ and
$\fW=\Bl_{\Xi}\fU$ is the blowup of $\fU$ in $\Xi$ and $W_0$ is the proper transform of $U_0$ in $\fW$ and $\fc{\fW} {W_0}$ is the formal completion of $\fW$ along $W_0$.
Since we have a natural map $\fc \fU {\Xi}\to \fc \fU {\Xi'}$, we may assume $\Xi=\Xi'$. 

\begin{claim}
We may assume that $U_0$ is irreducible.
\end{claim}

By Lemma \ref{lem;DecompositionXi}, we may assume that 
$\fUI=\underset{1\leq i\leq r}{\sqcup} U_i$,
where $U_i\subset \fUI$ are clopen and irreducible. 
Then, letting $\fU_i=\fc \fU {U_i}$, 
\[ \underset{1\leq i\leq r}{\sqcup} (\fU_i/\fY) \to (\fU/\fY)\]
is a Ninevich covering. Hence, the claim follows from Lemma \ref{lemm:mNisNormalForm}.

By the claim, there are two cases to be considered.


\def\th{\tilde{h}}

Case 1: $|\fUI|= |\Xi|$.
In this case, we have $\fc \fU {\Xi}=\fc \fU {\fUI}=\fU$.
So the trivial covering 
$\{(\fU/\fY) \to (\fU/\fY)\}$ gives a square \eqref{equa:VZUVWU}.

Case 2: $\fUI\not=\varnothing$ and $\Xi$ is nowhere dense in $\fUI$.
We will build a square \eqref{equa:VZUVWU} for some $\fY'\to \fY$ in $\Mod(\fY)$ with \eqref{equa:VZUVWU0} a trivial covering.
Let $\Sigma\subset Y_0=\uSpec(\cO_\fY/\cI)$ be the closure of the image of $\Xi$ under $U_0\to Y_0$. Since $\Xi$ is nowhere dense in $\fUI$ and $\fUI\to \fYI$ is \'etale, $\Sigma$ is a nowhere dense in $\fYI$. 

\begin{claim}\label{claim1;prop:cocont}
We may assume  $\Sigma\times_\fY\fU=\Xi\sqcup \Theta$ for a clopen $\Theta \subset \Sigma\times_\fY\fU$.
\end{claim}
\begin{proof}
By Lemma \ref{lem2;DecompositionXi}, there is a nowhere dense closed subscheme $D$ in $\Sigma$ such that letting $\fY'=\Bl_D(\fY)\in \Mod(\fY)$ and $\fU'=\fY'\times_\fY \fU$ and $\Sigma'$ 
(resp. $\Xi'$) be the proper transform of $\Sigma$ in $\fY'$ (resp. $\Xi$ in $\fU'$),
we have a decomposition $\Sigma'\times_{\fY'} \fU'=\Xi'\sqcup \Theta'$ for a clopen subscheme $\Theta'$. Then, it suffices to show that $\Bl_{\Xi'}(\fU')\to \fU'\to \fU$
factors through $\Bl_\Xi(\fU)$.
Note $\fU'=\Bl_D(\fY)\times_\fY \fU=\Bl_{\pi^{-1}(D)}(\fU)$ with the projection $\pi:\fU\to \fY$. 
Letting $h:\fU'\to \fU$ be the natural map and $E=h^{-1}(\pi^{-1}(D))$ be the exceptional divisor of $h$, the ideal sheaf $\cI_{\Xi'}$ is equal to $\cI_{h^{-1}(\Xi)}\cdot (\cI_E)^{-1}$.
So, we get the desired map 
$\Bl_{\Xi'}(\fU')=\Bl_{h^{-1}(\Xi)}(\fU')\to \Bl_\Xi(\fU)$ by the universality of blowups.
\end{proof}

Let $f:\tfY=\Bl_\Sigma(\fY) \to \fY$ and $\tfU=\tfY\times_\fY \fU$.
Using Claim \ref{claim1;prop:cocont}, we get a map 
\[\phi: \tfU=\Bl_\Sigma(\fY)\times_\fY \fU=
\Bl_{\Sigma\times_\fY\fU}(\fU)=\Bl_{\Xi\sqcup \Theta}(\fU)=
\Bl_{\Theta}(\Bl_{\Xi}(\fU))\to \Bl_\Xi(\fU), \]
 which factors through the projection $\tfU=\tfY\times_\fY \fU\to \fU$.
Let $\tY_0\subset \tfY$ be the proper transform of $Y_0$ in $\tfY$, Definition \ref{strict-peoprTransform}. Then, $\tU_0:=\tY_0\times_{\fY} \fU$ is the proper transform of $U_0$ in $\tfU=\Bl_{\Xi\sqcup \Theta}(\fU)$ so that $\phi(\tU_0)\subset W_0$, where $W_0$ is the proper transform of $U_0$ in $ \fW=\Bl_\Xi(\fU)$. Hence, $\phi$ induces a map
\[\fc {\tfY}{\tY_0}\times_\fY \fU =\fc {\tfU} {\tU_0} \to \fc \fW {W_0},\]
where $\fc {\tfY}{\tY_0}$ (rep. $\fc {\tfU} {\tU_0} $) is the formal completion of $\tfY$ along $\tY_0$ (resp. $\fU$ along $\tU_0$).
Thus, $\fY'=\fc {\tfY}{\tY_0}\in \Mod(\fY)$ satisfies the desired condition as in the beginning of Case 2. 
\end{proof}

\begin{coro} \label{coro:RZColimitLim}
Let $S$ be a Noetherian scheme. Then, the functor \eqref{rhoX}
\[ 
\rho^*_\fX: \Shv_{\procdh}(\fSchfp_S) 
\stackrel{}{\to}
\Shv(\RZ(\fX_{\Nis}))
\]
preserves colimits and finite limits, \cite[III.2.3]{SGA41}, \cite[00XL]{stacks-project}.
\end{coro}

Here are some lemmas that were used above.

\begin{lemm}\label{lem;DecompositionXi}
Let $\fY$ be a locally Noetherian formal scheme and $\pi:\fU\to \fY$ be an \'etale morphism. 
Assume that $|\fU|$ has only finitely many irreducible components.
Then, there exists a morphism $\fY'\to \fY$ in $\Mod(\fY)$ such that 
there exists a decomposition in $\fSch$
\[\fY'\times_\fY \fU=\underset{1\leq i\leq r}{\sqcup} \fU'_i\]
such that $|\fU'_i|$ is irreducible.
\end{lemm}
\begin{proof}
Let $Y=\uSpec(\cO_\fY/\cI)$ for an ideal $\cI$ of definition of $\fY$ and $U=Y\times_\fY\fU$.
By the assumption, $U$ has only finitely many irreducible components $U_1,\dots,U_r$.
Applying \cite[Lem.4.17]{KS23}, there exists a nowhere dense closed subscheme $D\subset Y$ such that
letting $Y'=\Bl_D(Y)$ and $U'=Y'\times_Y U$, there exist a decomposition 
$U'=\underset{1\leq i\leq r}{\sqcup} U'_i$, where $U'_i$ are clopen in $U'$.
Let $\tfY\to \fY$ be the blowup of $\fY$ in $D$, Definition \ref{defi;Blowup}.
Note that $Y'$ is the proper transform of $Y$ in $\tfY$, Definition \ref{strict-peoprTransform}. 
Let $\fY'$ be the formal completion of $\tfY$ along $Y'$.
Then, $\fY'\to \fY$ is in $\Mod(\fY)$ and $\fY'\times_{\fY}\fU$ is the formal completion of $\tfY\times_\fY\fU$ along $Y'\times_\fY \fU$ by \eqref{eq;formalCompletion}. Since 
\[Y'\times_\fY\fU=Y'\times_Y U=\underset{1\leq i\leq r}{\sqcup} U'_i,\]
this implies that $\fY'\times_{\fY}\fU$ is a disjoint sum of the formal completions of $\tfY\times_\fY\fU$ along 
$U'_i$, which completes the proof.
\end{proof}

\begin{lemm}\label{lem2;DecompositionXi}
Let $\fY$ be a locally Noetherian formal scheme and $\pi:\fU\to \fY$ be an \'etale morphism. Let $Y=\uSpec(\cO_\fY/\cI)$ for an ideal $\cI$ of definition of $\fY$ and $U=Y\times_\fY\fU$.
Let $\Xi\subset U$ be a closed subscheme and  $\Sigma$ be the closure of the image of $\Xi$ in $Y$. Then, there exists a nowhere dense closed subscheme $D$ in $\Sigma$ such that the blowup $\fY'=\Bl_D(\fY)$ of $\fY$ in $D$, Definition \ref{defi;Blowup}, satisfies the following condition: Let $\fU'=\fY'\times_\fY \fU$ and $\Sigma'$ 
(resp. $\Xi'$) be the proper transform of $\Sigma$ in $\fY'$ (resp. $\Xi$ in $\fU'$). 
Then, we have a decomposition $\Sigma'\times_{\fY'} \fU'=\Xi'\sqcup \Theta'$ for a clopen subscheme $\Theta'$.
\end{lemm}
\begin{proof}
Since $f:U\to Y$ is \'etale, we have
$f^{-1}(\Sigma)=\Xi\cup \Theta$ for a closed subscheme $\Theta$ and
$\Xi\sqcup \Theta \to \Sigma$ is generically \'etale.
By Raynaud-Gruson platification, \cite[081R]{stacks-project}, we find a blowup 
$\Sigma'=\Bl_D(\Sigma) $ of $\Sigma$ in a nowhere dense closed subscheme $D\subset \Sigma$ such that the proper transform $\Xi' \sqcup \Theta' \to \Sigma'$ of $\Xi \sqcup \Theta \to \Sigma$ is flat. %
By \cite[Lem. 4.15]{KS23} this implies $\Xi' \sqcup\Theta' \to \Sigma' \times_Y U$ is also flat. It is proper and an isomorphism over a dense open subset by construction, so it is in fact an isomorphism by \cite[Lem. 4.16]{KS23}. 
Let $\fY'=\Bl_D(\fY)$ be the blowup of $\fY$ in $D$ and $\fU'=\fY'\times_\fY \fU$.
By the construction, $\Sigma'\times_{\fY'} \fU'=\Sigma' \times_Y U\simeq \Xi'\sqcup \Theta'$, which completes the proof.
\end{proof}

\subsection{Finiteness of homotopy dimension}

In this subsection, we complete the proof of Theorem \ref{thm2;Hypercomplete}.

\begin{prop}[{cf.\cite[Prop.2.4.3]{EHIK}}] \label{prop:RZfinDim}
Let $\fX$ be a locally Noetherian formal scheme with finite Krull dimension $d$.
Then the ∞-topos $\Shv(\RZ(\fX),\Spc)$ has homotopy dimension $\leq d$.
\end{prop}
\begin{proof}
By \eqref{ShvRZ}, we have 
\[ \Shv(\RZ(\fX_{\Nis}),\Spc) = \varinjlim_{\fY\in \Mod(\fX)} \Shv(\fY_{\Nis},\Spc).\]
By \eqref{Nisequiv} and \cite[3.18]{CM21}, 
$\Shv(\fY_{\Nis,\Spc})$ has homotopy dimension $\leq d_\fY$, where $d_\fY$ is the Krull dimension of $\fY$.
By Remark \ref{rema:genXY}, the system of those $\fY\in \Mod(\fX)$ with $d_\fY\leq d$ is cofinal. Hence, the assertion follows from \cite[Cor.3.11]{CM21},
\def\bir{\mathrm{bir}}
\end{proof}

We now start proving Theorem \ref{thm2;Hypercomplete}.
We follow the proof of \cite[Th.2.4.15]{EHIK}.
The proof proceeds by induction on $d$. Let $F\in \Shvzen(\fSchffpSX,\Spc)$ be a $d$-connective object. We want to show that $F(\fX)$ is nonempty. 
If $d < 0$, then $\fX$ is empty and the assertion is trivial. 
Assume $d\geq 0$. By Corollary \ref{coro:RZColimitLim}, the functor \eqref{rho}
\[ \rho^*: \PSh(\fSchffpSX,\Spc) \to \PSh(\RZ(\fX_{\Nis}),\Spc)\]
induced by $\rho: \RZ(\fX_{\Nis}) \to \fSchffpSX; (\fU/\fY) \mapsto \fU$,
preserves colimits and finite limits. Hence, by \cite[Prop. 6.5.1.16(4)]{HTT} (see also \cite[\S6.2]{KS23}), $\rho^*$ sends $d$-connective objects to $d$-connective objects.
Since $\Shv(\RZ(\fX_{\Nis}),\Spc)$ has homotopy dimension $\leq d$ by
Proposition \ref{prop:RZfinDim}, $\rho^*F$ has a global section.
We have 
\[(\rho^*F)(\fX/\fX) =a(F \circ \rho)(\fX/\fX) = \varinjlim_{\fY \in \Mod(\fX)} F(\fY),\]
where $a: \PSh(\RZ(\fX_{\Nis}),\Spc)\to \Shv(\RZ(\fX_{\Nis}),\Spc)$ is the sheafification functor and the second equality follows from Lemma \ref{prop:sheafificationCalc}.
Hence, there is $f:\fY\to \fX \in \Mod(\fX)$ such that $F(\fY)\not=\varnothing$.
Let $i:Z\subset \fX$ be a nowhere-dense closed subscheme such that
$f$ belongs to $\Mod_{\fX\setminus Z}(\fX)$, Definition \ref{def;Modificcation}.
Note $\dim(Z)<d$. By Remark \ref{rema:genXY}, we may assume $\dim(E)<d$ for $E=Z\times_\fX \fY$. By Lemma \ref{excisionAdmMod}, we have 
\[ F(\fX) \cong F(\fY) \times_{F(\fc \fY E)} F(\fc \fX Z). \]
By the induction hypothesis, ${\fSchffpS}_{/\fM}$ with $\fM=\fc \fX Z$ or $\fc \fY E$ has homotopy dimension $<d$. Since the restrictions of $F$ to ${\fSchffpS}_{/\fM}$ remains $d$-connective by \cite[Prop. 6.5.1.16(4)]{HTT} (see also \cite[Example 6.4]{KS23}), we conclude that  $F(\fc \fX Z)$ and $F(\fc{\fY} E)$ are connected
by \cite[Def.7.2.1.6, Lem.7.2.1.7]{HTT} (see also \cite[Rem. 7.6]{KS23}). 
Hence $F(\fX)$ is non-empty and the proof is complete.

\begin{coro} \label{coro:finCohDim}
Let $\fX$ be a Noetherian formal scheme of Krull dimension $d \geq 0$.
For any sheaf of abelian groups $F \in \Shv_{\procdh}(\fSchffp_\fX, \Ab)$ we have
\[ H_{\procdh}^n(\fX, F) = 0; \qquad n > d. \]
\end{coro}

\begin{proof}
This is \cite[Cor.7.2.2.30]{HTT}. Cf. also \cite[Def.7.2.2.14, Rem.7.2.2.17]{HTT}.
\end{proof}

\section{Fibre functors}

In this section, we fix a Noetherian adic ring $S$
and characterise fibre functors of the 1-topos $\Shvzen(\fSchffpS)=\Shvzen(\fSchffpS.\Set)$. First, we recall the following.

\begin{defi} \label{defi:FF}
Let $(C,\tau)$ be a site admitting finite limits.
Recall that a \emph{fibre functor} of an $1$-topos $\Shv_τ(C)$ of sheaves of sets, is a continuous morphism of topoi $\phi^*: \Shv_τ(C) \rightleftarrows \Set: \phi_*$, or equivalently, a functor $\phi^*: \Shv_τ(C) \to \Set$ which preserves colimits and finite limits. Let $\Fib(\Shv_τ(C))$ denote the category of fiber functors of $\Shv_τ(C)$.
\end{defi}

\def\aff{\mathrm{aff}}

Let $\AffSchffpS\subset \fSchffpS$ be the full subcategory of affine formal schemes $\Spf(A)$ for $A\in \ARingffp_S$, Definition \ref{def;catftadicring}.
Let $τ$ be a topology on $\fSchffpS$ which is finer than the Zariski topology.
Then, the canonical functor $\Shv_\tau(\fSchffpS) \to \Shv_{\tau^{\aff}}(\AffSchffpS)$ is an equivalences, where $\tau^{\aff}$ is the topology on $\AffSchffpS$ induced by $\tau$, i.e. the finest topology such that the image in $\fSchffpS$ of any $\tau^{\aff}$-covering family is a $\tau$-covering family.
So, there is an equivalence of categories of fiber functors
\eq{eqFib}{
\Fib(\Shv_\tau(\fSchffpS)) \simeq \Fib(\Shv_{\tau^{\aff}}(\AffSchffpS)).}
Hence, by \cite[Pro.7.13]{Joh77}, there is a bijection between fibre functors of $\Shv_\tau(\fSchffpS)$ and cofiltered pro-objects 
\eq{Pbullet}{P_\bullet=\catprojlim {\lambda\in \Lambda^{op}} P_\lambda \qwith 
P_\lambda=\Spf(R_\lambda)\in \AffSchffpS}
indexed by the opposite category of a filtered category $\Lambda$, which satisfies the $\tau$-locality condition: For every $τ$-covering 
$\{\fY_i \rmapo{u_i} \fX\}_{i \in I}$, the morphism of sets
\begin{equation}\label{homPbullet}
\coprod_{i \in I} \varinjlim_{\lambda\in \Lambda} \Hom_{\fSchS}(P_\lambda, \fY_i) \to \varinjlim_{\lambda\in \Lambda} \Hom_{\fSchS}(P_\lambda,  \fX)
\end{equation}
is surjective. Let 
\eq{Rbullet}{R_\bullet=\catinjlim {\lambda\in \Lambda}  \Rl \in \Ind(\ARingffp_S)}
be the corresponding ind-objects and $(\Rinf,\fn)=\Psi(R_\bullet)$, Definition \ref{def;Phi}.
The fibre functor associated to a pro-object $P_\bullet$ \eqref{Pbullet} is given by
\[ \phi_{R_\bullet}: \PSh(\fSchffp_S) \to \Set\;;\; F \to 
\varinjlim_{\lambda\in \Lambda} F(P_\lambda),\]
which factors through $ \Shv_{\tau}(\fSchffp_S) $ by the $\tau$-locality.



\begin{prop}\label{prop;fibrePsi}
Let $R_\bullet, R'_\bullet \in \Ind(\ARingffp_S)$.
If $\Psi(R_\bullet)=\Psi(R'_\bullet)$, there is a natural equivalence
$\phi_{R_\bullet}\simeq \phi_{R'_\bullet}$ of functors. 
\end{prop}
\begin{proof}
We want to define a bijection $\xi_F: \phi_{R_\bullet}(F)\simeq \phi_{R'_\bullet}(F)$ functorial in $F\in \PSh(\fSchffp_S) $. Since $\phiR$ commutes with colimits, it is enough to construct $\xi_F$ for representable objects $F=\Hom_{\ARingffp_S}(A,-)$ for $A\in\ARingffp_S$. This follows from Corollary \ref{cor;lem;reduction}.
\end{proof}

By the above proposition, Theorem \ref{prop-intro:procdhLocal} follows from the following.

\begin{theo}[Characterisation of fibre functors] \label{prop:procdhLocal}
A cofiltered pro-object $P_\bullet$ \eqref{Pbullet} is pro-cdh local if and only if $\Rinf/\fn$ is a henselian valuation ring and
\eq{eq;prop:procdhLocal}{\Rinf\simeq \Rinf/\fn\times_{\kn} \Rn,}
where $\kn$ is the residue field of $\fn$ and $\Rn$ is the localization of $\Rinf$ at $\fn$.
In other words, $\Rinf \subseteq \Rn$ is the set of elements whose residue classes mod $\fn$ lie in $\Rinf/\fn$.
\end{theo}

\begin{rema}\label{rema0;prop:procdhLocal}
Since $(\Rinf,\fn)$ is a henselian pair by Lemma \ref{lem;Phi}, the condition of the theorem implies that $\Rinf$ is henselian local.
\end{rema}

\begin{rema}\label{rema;prop:procdhLocal}
By Deligne's completeness theorem, \cite[Prop.VI.9.0]{SGA42} or \cite[Thm.7.44, 7.17]{Joh77}, those fibre functor of $\Shv_{\procdh}(\fSchffp_S)$ satisfying the condition of Theorem \ref{prop:procdhLocal} form a conservative family, i.e. a morphism $f$ in $\Shvzen(\fSchffpS)$ is an isomorphism if and only if $\phi_{P_\bullet}(f)$ is an isomorphism for all such $\phi_{P_\bullet}$.
Equivalently, a family $\{\fY_i \to \fX\}_{i \in I}$ in $\fSchffp_S$ is a pro-cdh covering if and only if 
$\sqcup_{i \in I} \phi_{P_\bullet}(\fY_i) \to \phi_{P_\bullet}(\fX)$ is surjective for all such $\phi_{P_\bullet}$, \cite[Exposé IV, Prop.6.5(a)]{SGA41}, 

\end{rema}

We need a preliminary for the proof of Theorem \ref{prop:procdhLocal}.
Let 
\[\hAS 2=\Spf(\cO_S\{x,y\}\})=\AZ 2 \times_{\Spec(\ZZ)} S\]
with $\cO_S\{x,y\}=\varprojlim_{m} \cO_S/I^m[x,y]$ 
and $\AZ 2=\Spec \ZZ[x,y]$ and the fiber product is taken in $\fSch$ (cf. Definition \eqref{hotimes}).
Let $\hASo 2$ be the formal completion of $\AS 2$ along the origin $0=\{x=y=0\}$, Definition \ref{def;formalCompletion}.  
Note $\hASo 2=\Spf(\cO_S[[x,y]])$, where 
\[\cO_S[[x,y]]=\varprojlim_{m,l} \cO_S/I^m[x,y]/(x,y)^l,\]
is the $I\cO_S\{x,y\}+(x,y)$-adic completion of $\cO_S[x,y]$.
Let $\BlAt$ is the blowup of $\AZ 2$ at $0$ and $\hBlAtS=\BlAtZ\times_{\AZ 2} \hAS 2$.


\begin{lemm}\label{lem1:procdhLocal}
Let $\tau$ be the topology on $\fSchffpS$ generated by Nisnevich topology and the following covering families:
\begin{enumerate}
\item[(SmB)] Smooth blow up covering\;
\[ \left \{ \hBlAtS \rightarrow \hAS 2 , \; \hASo 2 \rightarrow \hAS 2 \right \}. \]

\item[(Axe)] Letting $\cA=\cO_S\{x,y\}/(xy)$, 
\[ \left \{ \Spf(\cA/(y)) \rightarrow \Spf(\cA), \; \Spf(\cAx)  \rightarrow \Spf(\cA) \right \}\]
where $\cAx$ is the $(x)$-adic completion of $\cA$.
Note $\cAx=\cO_S\{y\}[[x]]/(xy)$, where $\cO_S\{y\}[[x]]$ is the $I\cO_S\{x,y\}+(x)$-adic completion of $\cO_S\{x,y\}$.
\end{enumerate}
A cofiltered pro-object \eqref{Pbullet} is $\tau$-local if and only if $\Rinf$ satisfies the condition of Theorem \ref{prop:procdhLocal}.
\end{lemm}

\begin{rema}\label{rem;generatorTopology}
Clearly, the Nisnevich topology and (SmB) generate the same topology as the Nisnevich topology and (SmB') do, where 
\begin{enumerate}
\item[(SmB')] \;
\[ \left \{ \Spf(\cO_S\{x,y/x\}) \sqcup \Spf(\cO_S\{y,x/y\})  \rightarrow \hAS 2 , \;
  \hASo 2 \rightarrow \hAS 2  \right \}. \]
\end{enumerate}
\end{rema}

\medbreak\noindent
{\it Proof of Lemma \ref{lem1:procdhLocal}:}
Using Remark \ref{rem2;lem;Phi}, we compute
\[ \varinjlim_\lambda \Hom_{\fSchS}(P_\lambda,  \Spf(A))=
\varinjlim_\lambda\Hom_{\pARing_S}(A,\Rl)\]
for $A=S\{x,y\}$, $S[[x,y]]$, $\cA=S\{x,y\}/(xy)$ and $\cAx=S\{y\}[[x]]/(xy)$ as follows: 
\eq{SxyRbullet}{
\varinjlim_\lambda \Hom_{\fSchS}(P_\lambda,  \Spf(S\{x,y\})) \simeq \Rinf\times \Rinf.}
\eq{SxyxyRbullet}{
\varinjlim_\lambda \Hom_{\fSchS}(P_\lambda,  \Spf(S[[x,y]])) \simeq \fn\times \fn,}
\eq{SxyRbulletxy}{
\varinjlim_\lambda \Hom_{\fSchS}(P_\lambda,  \Spf(\cA)) \simeq 
\{(a,b)\in \Rinf\times \Rinf|\; ab=0\},}
\eq{SxyxRbulletxy}{
\varinjlim_\lambda \Hom_{\fSchS}(P_\lambda,  \Spf(\cAx)) \simeq 
\{(a,b)\in \fn\times \Rinf|\; ab=0\}.}
In view of \eqref{SxyRbullet} and \eqref{SxyxyRbullet}, the surjectivity of \eqref{homPbullet} for the covering (SmB') in Remark \ref{rem;generatorTopology} is equivalent to that
\begin{enumerate}
 \item[{$(*)$}] $\forall a, b \in \Rinf$, we have $a | b$ or $b|a$ or $a$ and $b$ are both in $\fn$.
\end{enumerate}
In view of \eqref{SxyRbulletxy} and \eqref{SxyxRbulletxy},
the surjectivity of \eqref{homPbullet} for the covering (Axe) is equivalent to the condition
\begin{enumerate}
 \item[{$(**)$}] all zero divisors of $\Rinf$ are in $\fn$.
\end{enumerate}
$(*)$ implies that $\forall a, b \in \Rinf/\fn$, we have $a | b$ or $b|a$ so that  $\Rinf/\fn$ is a valuation ring. $(**)$ implies that $\Rinf\to \Rn$ is injective, proving the injectivity of \eqref{eq;prop:procdhLocal}. To show its surjectivity, consider the commutative diagram
\[\xymatrix{
\fn \ar[r]\ar[d] & \Rinf \ar[r]\ar[d] & \Rinf/\fn \ar[d] \\
\fn\Rn \ar[r] & \Rn \ar[r] &\kn\\}\]
where the rows are exact. So, it suffices to show the surjectivity of 
$\fn\to \fn\Rn$. Take $a/r\in \fn\Rn$ with $a\in \fn$ and $r\in \Rinf\setminus\fn$. By $(*)$, either $a|r$ or $r|a$.
In the former case, we have $r\in \fn$ contradicting $r\not\in \fn$.
In the latter case, we can write $a=qr\in \fn$ for some $q\in \Rinf$ and then $q\in \fn$ since $\fn$ is prime and $r\not\in \fn$. Since $q$ maps to $a/r$ under $\fn\to \fn\Rn$, this proves the desired surjectivity, proving the isomorphism\eqref{eq;prop:procdhLocal}. Conversely, it is easy to check that \eqref{eq;prop:procdhLocal} implies the conditions $(**)$ and $(**)$.

Next, we prove that $\Rinf/\fn$ is henselian local. It suffices to prove that 
$\Rinf$ is henselian local. 
By what we have shown above, $\Rinf$ is local. By \cite[10.153.1]{stacks-project}, it suffices to prove the following.

\begin{claim}
Let $\fm\subset\Rinf$ be the maximal ideal. Assume given a monic $f\in \Rinf[t]$ and $\alphab\in \Rinf/\fm$ such that $\fb(\alphab)=0$ and $\fb'(\alphab)\not=0$, where 
$\fb$ is the image of $f$ in $\Rinf/\fm[t]$. Then, there exists $\alpha\in \Rinf$ such that 
$f(\alpha)=0$ and $\alpha$ maps to $\alphab$.
\end{claim}

For each $\lambda\in \Lambda$, let $\fm_\lambda$ be the inverse image of $\fm$ in $R_\lambda$. The assumption implies that there exists $\lambda$ such that $f$ (resp. $\alphab$) is the image of $\fl\in \Rl[t]$ (resp. $\abl\in \Rl/\fml$) such that $\fbl(\abl)=0$ and 
$\fbl'(\abl)\in (\Rl/\fm_\lambda)^\times$.
Letting $B_\lambda=\Rl[t]/(\fl)$, we have a decomposition 
$B_\lambda/\fml\simeq \Rl/\fml\times Q$ as a ring, where the projection $\pi:B_\lambda/\fml=\Rl/\fml[t]/(\fbl)\to \Rl/\fml$ is given by the substitution of $\abl$.
By \cite[10.144.2(1)]{stacks-project}, there exists $h\in B_\lambda\setminus \Ker(\pi)$ such that $\xi: \Rl\to C:=B_\lambda[h^{-1}]$ is \'etale and $\xi\otimes \Rl/\fml$ is an isomorphism.
Letting $\hC$ be the $\fn_\lambda C$-adic completion of $C$,
\[\hV\sqcup \big(\Spf(\Rl)\setminus \Spec(\Rl/\fml)\big)\to \Spf(\Rl)\qwith \hV=\Spf(\hC)\]
is a Nisnevich covering from Definition \ref{def;fXnis}. 
Let $\Laml\subset \Lambda$ be the subcategory of objects under $\lambda$. Noting that $\Laml$ is cofinal in $\Lambda$,
the surjectivity of \eqref{homPbullet} for this covering implies that the canonical element of 
$\varinjlim_{\alpha\in \Laml} \Hom(\Rl,R_\alpha)$ lifts to 
$\varinjlim_{\alpha\in \Laml}\Hom(\hC,R_\alpha)$, which implies that the natural map $\iota_\lambda:\Rl\to \Rinf$ lifts to a map $g: \hC \to \Rinf$.
Composing $g$ with $B_\lambda\to C \to \hC$, we get a map $\psi: B_\lambda\to \Rinf$ whose composite with $\Rl\to B_\lambda$ is $\iota$. Then, $\alpha=\psi(t)\in \Rinf$ satisfies the desired condition of the claim.

To complete the proof of Lemma \ref{lem1:procdhLocal}, it remains to show that for a cofiltered pro-object \eqref{Pbullet} satisfying the condition of Theorems \ref{prop:procdhLocal}, \eqref{homPbullet} is surjective for the Nisnevich coverings. 
By Proposition \ref{prop;NisfX}, it suffices to treat the Nisnevich covering $\fU\sqcup \fY \to \fX$, where $\fX=\Spf(A)$ and $\fY=\Spf(B)$ are affine, $j:\fU\to \fX$ is an open immersion and $f:\fY\to \fX$ is \'etale such that 
writing $X_0=\Spec(A/I)$ and $Y_0=\Spec(B/IB)$ for an ideal $I\subset A$ of definition, the induced morphism $Y_0\to X_0$ of schemes is an isomorphism over a closed subscheme $Z\subset X_0$ such that $|Z|=|\fX|\setminus |\fU|$. 
Letting $U=\Spec(A)\setminus Z$, $\fU$ is the formal completion of $U$ along $U_0=U\times_{\Spec(A)} X_0$.
By Theorem \ref{thm;algebrizaitonEtale}, we may further assume that there is an \'etale $A$-algebra $B'$ whose $IB'$-adic completion $\hat{B'}$ is isomorphic to $B$.

Choose $\lambda\in \Lambda$ and $\phi\in \Hom_{\fSchS}(P_\lambda,  \fX)=\Hom_{\ARing_S}(A,\Rl)$, and let 
\eq{tildephi}{\tilde{\phi} \in \varinjlim_{\alpha\in \Laml} \Hom_{\ARing_S}(A,R_\alpha)=
\varinjlim_{\alpha\in \Laml}\Hom_{\fSchS}(P_\alpha,  \fX)}
be the image of $\phi$. It suffices to show $\tilde{\phi}$ lifts to 
$\varinjlim_{\alpha\in \Laml}\Hom_{\fSchS}(P_\alpha, \fU\sqcup \fY)$.
Let $\phi_\infty\in \Hom_{\pARing_S}(A,\Rinf)$ be the image of $\phi$. 
Let $\phi_\infty$ also denote the induced map $\Spec(\Rinf) \to \Spec(A)$.

Assume that the image of $\phi_\infty$ is contained in $U$.
Since $\Rinf$ is local, $\phi_\infty$ factors through an affine open
$\Spec(C)$ of $U$, where $C$ is an algebra of finite type over $A$.
Let $\psi:C\to \Rinf$ be the corresponding map of rings which lifts $\phi_\infty$.
By Lemma \ref{lem1;adicring}, we have $\psi\in  \Hom_{\pARing_A}(C,\Rinf)$, where $C$ is equipped with the $IC$-adic topology. 
By Proposition \ref{lem;compactness}, this means that $\tilde{\phi}$ lifts to
\[\varinjlim_{\alpha\in \Laml} \Hom_{\ARing_S}(\hC,R_\alpha)
= \varinjlim_{\alpha\in \Laml}\Hom_{\fSchS}(P_\alpha,  \Spf(\hC)).\]
Since $\Spf(\hC)$ is an affine formal open of $\fU$ by the construction, we conclude that $\tilde{\phi}$ lifts to
$\varinjlim_\alpha\Hom_{\fSchS}(P_\alpha, \fU)$.

If the image of $\phi_\infty$ is not contained in $U$, then the image of 
$\phib_\infty: \Spec(\Rinf/\fm)\to \Spec(\Rinf) \to \Spec(A)$ is contained in $Z\subset X_0$, where $\fm$ is the maximal ideal of $\Rinf$.
So, $\phib_\infty$ lifts to a map $\psib: \Spec(\Rinf/\fm)\to Y_0\to \Spec(B')$.
Since $\Rinf$ is henselian (cf. Remark \ref{rema0;prop:procdhLocal}) and $A\to B'$ is \'etale, $\psib$ lifts to a map $\psi: B' \to \Rinf$ which lifts $\phi_\infty$. 
By Lemma \ref{lem1;adicring}, we have $\psi\in  \Hom_{\pARing_A}(B',\Rinf)$, where $B'$ is equipped with the $IB'$-adic topology. 
By Proposition \ref{lem;compactness}, this means that $\tilde{\phi}$ lifts to
\[\varinjlim_{\alpha\in \Laml} \Hom_{\ARing_S}(B,R_\alpha)
= \varinjlim_{\alpha\in \Laml}\Hom_{\fSchS}(P_\alpha,  \fY).\]
This completes the proof of Lemma \ref{lem1:procdhLocal}. 
\medbreak

\def\tpsi{\tilde{\psi}}

Now, we prove Theorem \ref{prop:procdhLocal}. It suffices to prove that every pro-cdh covering is a 
$\tau$-covering. By Lemma \ref{lem1:procdhLocal} and Remark \ref{rema;prop:procdhLocal} (applied to the $\tau$-topology), it is enough to show that for a cofiltered pro-object \eqref{Pbullet} satisfying the condition of Theorems \ref{prop:procdhLocal}, \eqref{homPbullet} is surjective for  the covering $\fc \fX Z\sqcup \fY \to \fX$ from (FAB) in Definition \ref{defi:procdh}. 
By \eqref{eqFib}, we may assume that $\fX=\Spf(A)$ is affine and $Z\subset \fX$ is defined by an ideal $K\subset A$ containing an ideal $I$ of definition of $A$. 
Moreover, by Definition \ref{defi:locallyalgebriaable}, we may further assume that 
there is a morphism $g:Y\to X=\Spec(A)$ of schemes which is proper of finite type such that $\fY=Y\times_X \fX$ and that $g$ is an isomorphism over $\Spec(A)\setminus V(J)$ for some ideal $J\subset A$ such that $K\subset J+ \sqrt{I}$.

Choose $\lambda,\phi,\tilde{\phi}$ as in \eqref{tildephi}.
 It suffices to show that $\tilde{\phi}$ lifts to 
 \[\varinjlim_{\alpha\in \Laml}\Hom_{\fSchS}(P_\alpha, \fc \fX Z\sqcup \fY ).\]
Let $\phi_\infty\in \Hom_{\pARing_S}(A,\Rinf)$ denote the image of $\phi$. 
By \eqref{eq;prop:procdhLocal}, we have $\Rinf\simeq \Rinf/\fn\times_{\kn} \Rn$.
Let $\phi_1: A\to \Rinf/\fn$ and $\phi_2: A\to \Rn$ be the composite maps of $\phi_\infty$ and the projections.

If $\phi_\infty(J)\not\subset \fn$, $\phi_2$ factors through $A[h^{-1}] \to\Rn$ for some $h\in J$ since $\Rn$ is local.
Thus, the induced map $\Spec(\Rn) \to \Spec(A)$ factors through $\Spec(A)\setminus V(J)$ so that it lifts to $\psi_2:\Spec(\Rn) \to Y$. 
Noting that $\Rinf/\fn$ is a valuation ring of $\kn$ and $Y\to \Spec(A)$ is proper, the composite
$\Spec(\kn)\to \Spec(\Rn) \to Y$ factors through $\psi_1:\Spec(\Rinf/\fn) \to Y$ which lifts the map $\Spec(\Rinf/\fn) \to \Spec(A)$ induced by $\phi_1$. By \eqref{eq;prop:procdhLocal}, we have
\begin{multline*}
 \Hom_{\Sch_S}(\Spec(\Rinf),-)=\Hom_{\Sch_S}(\Spec(\Rinf/\fn,-)\sqcup_{\Spec(\kn)}\Spec(\Rn),-)\\
=
\Hom_{\Sch_S}(\Spec(\Rinf/\fn),-)\times_{\Hom_{\Sch_S}(\Spec(\kn),-)}
\Hom_{\Sch_S}(\Spec(\Rn),-).\end{multline*}
Thus, $\psi_1$ and $\psi_2$ give rise to a map $\psi:\Spec(\Rinf)\to Y$ which lifts the map $\Spec(\Rinf)\to \Spec(A)$ induced by $\phi_\infty$.
Since $\Rinf$ is local, $\psi$ factors through an affine open
$\Spec(B)$ of $Y$, where $B$ is an algebra of finite type over $A$.
Let $\chi: B\to \Rinf$ be the corresponding map of rings which lifts $\phi_\infty$.
By Lemma \ref{lem1;adicring}, we have $\chi\in  \Hom_{\pARing_S}(B,\Rinf)$, where $B$ is equipped with the $IB$-adic topology, and $\chi$ lifts $\phi_\infty$.
By Proposition \ref{lem;compactness}, this means that $\tilde{\phi}$ from \eqref{tildephi} lifts to
\[\varinjlim_{\alpha\in \Laml} \Hom_{\ARing_S}(\hB,R_\alpha)
= \varinjlim_{\alpha\in \Laml}\Hom_{\fSchS}(P_\alpha,  \Spf(\hB)).\]
where $\hB$ is the $IB$-adic completion $\hB$ of $B$.
Since $\Spf(\hB)$ is an affine formal open of $\fY$ by the construction, we conclude that $\tilde{\phi}$ lifts to
$\varinjlim_\alpha\Hom_{\fSchS}(P_\alpha, \fY)$.

\def\AK{A_{\{K\}}}

Next assume $\phi_\infty(J)\subset \fn$. By Lemma \ref{lem;Phi}, we have
$\phi_\infty(\sqrt{I})\subset \fn$ so that $\phi_\infty(K)\subset \phi_\infty(J+\sqrt{I})\subset \fn$.
Let $\AK \in \pARing_S$ be the ring $A$ equipped with $K$-adic topology.  
By Lemma \ref{lem;cartesian}, we have a cartesian square
\[\xymatrix{\varinjlim_{\alpha\in \Laml} \Hom_{\pARing_S}(\AK,R_\alpha)\ar[r]\ar[d]
& \Hom_{\pARing_S}(\AK,\Rinf)\ar[d] \\
\varinjlim_{\alpha\in \Laml} \Hom_{\Ring_S}(A,R_\alpha)\ar[r] &\Hom_{\Ring_S}(A,\Rinf)\\}\]
By \eqref{eq1;rem;def;adicring}, we have a bijection
\[ \varinjlim_{\alpha\in \Laml}  \Hom_{\pARing_S}(\AK,R_\alpha)=
\varinjlim_{\alpha\in \Laml}  \Hom_{\ARing_S}(\fc A K,R_\alpha)=\varinjlim_{\alpha\in \Laml}  \Hom_{\fSchS}(P_\alpha,  \fc \fX Z).\]
Hence, $\phi_\infty(K)\subset \fn$ implies that $\tilde{\phi}$ lifts to $\varinjlim_\alpha \Hom_{\fSchS}(P_\alpha,  \fc \fX Z)$. 
This completes the proof of Theorems \ref{prop:procdhLocal}.

\begin{coro}\label{thm;generatorTopology}
The pro-cdh topology on $\fSchffpS$ is generated by the Nisnevich topology and the covering families (SmB) and (Axe) from Lemma \ref{lem1:procdhLocal}.
\end{coro}


\section{A topos-theoretic interpretation of Weibel's vanishing}\label{Weibel}
\def\tF{\tilde{F}}

We set $\cC = \Spt$ or $D(\ZZ)$. For an integer $N$, $\cC_{\geq N}$ denotes the full subcategory of objects of $\cC$ supported in homological degree $i\geq N$.  
Let $\Sch^\qcqs$ (resp.$\Sch^\qcqs_\FF$) be the category of qcqs schemes
(resp. qcqs schemes over a field $\FF$).

\begin{theo} \label{theo:procdhWebelVanishing}
For $F \in \PSh(\Sch^{\qcqs}, \cC)$ and integer $N$, consider the conditions.%
\begin{enumerate}
 \item[(Desc)] For every Noetherian scheme $X$, the restriction $F$ to the category of schemes of finite type over $X$ is a pro-cdh sheaf from Definition \ref{def-intro;procdhSch}.
 \item[(Fin)] $F$ is finitary, in the sense that it preserves filtered colimits of rings.
 \item[(Rig)$_N$] 
$\fib(F(A)\to F(A/I))\in \cC_{\geq -N}$ for every Noetherian ring $A$ and nilpotent ideal $I\subset A$.
\item[(Val)$_N$] $F(R)\in \cC_{\geq -N}$ for every henselian valuation ring $R$.
\end{enumerate}
If $F$ satisfies (Desc), (Fin), (Rig)$_N$, and (Val)$_N$,
then for every Noetherian scheme $X$ of finite Krull dimension $d$, we have
\begin{equation*} \label{eq;rmk;apcdhWebel}
\pi_i F(X)=0 \qfor  i<-N-d.
\end{equation*}
\end{theo}

We note that the essentially same statement as above is proved by Elmanto-Morrow \cite[Prop.8.10]{EM23} using an argument used in the proof of Weibel's conjecture \cite[Th.B]{KST-Weibel} and also an idea used in the proofs of Weibel and Vost's conjectures for finite type schemes over fields of characteristic zero 
\cite{CHSW08}, \cite{CHW08}.
We will give a new proof using the main results of this paper, Theorem~\ref{prop:procdhLocal} and Corollary \ref{coro:finCohDim}.

\begin{proof}
For $X$ as above, let $\fSchffpX$ be the category of formal schemes formally of finite type over $X$ (viewed as a discrete formal scheme).
Let $\Fc$ be the presheaf on $\fSchffpX$ defined by
\eq{eq2;continutity}{ \Fc(\fY) =\varprojlim_{\cI} F(\fY_\cI)\qfor \fY\in \fSchffpX,}
where $\cI$ ranges over ideals of definition of $\fY$. 
By \eqref{Shvprocdhcont}, (Desc) implies $\Fc\in \Shvzen(\fSchffpX)$. 
Hence, we have the descent spectral sequences
\eq{DescentSS}{ E_2^{p,q} = H^p_{\procdh}(X,a_{\procdh}\pi_{-q}\Fc) \Rightarrow \pi_{-p-q}F(X),}
where $a_{\procdh}\pi_{i}\Fc$ is the pro-cdh sheafification of the presheaf $\pi_i \Fc$ of abelian groups on $\fSchffpX$.
By Corollary \ref{coro:finCohDim}, we have $E_2^{p,q} =0$ for $p>d$.
So, it suffices to show that $a_{\procdh}\pi_{i}\Fc=0$ for $i<-N$.
By Theorem~\ref{prop:procdhLocal}, it is enough to show 
\[ \varinjlim_\lambda\pi_i\Fc(R_\lambda) =0\qfor i <-N\]
for every ind-system $\catinjlim \lambda R_\lambda$ of adic rings formally of finite type over $X$ such that $\Rinf=\varinjlim_{\lambda} R_\lambda$ satisfies the condition of the theorem.
Let $\fn_\lambda\subset R_\lambda$ be the ideal of topologically nilpotent elements and $\fn=\varinjlim_{\lambda} \fn_\lambda$. For $i<-N$, we have
\[ \begin{aligned}
\varinjlim_\lambda\pi_i\Fc(R_\lambda) &\overset{\eqref{eq2;continutity}}{=} \varinjlim_\lambda \pi_i\varprojlim_n F(R_\lambda/\fn_\lambda^n) =  \varinjlim_\lambda\pi_i F(R_\lambda/\fn_\lambda)=\pi_iF(\Rinf/\fn) =0\\
\end{aligned}
\]
where the second equality follows from (Rig)$_N$ in view of Milnor's exact sequence\footnote{For $i=-N-1$, (Rig)$_N$ implies $\pi_{i+1}F(R_\lambda/\fn_\lambda^{n+1})\to \pi_{i+1}F(R_\lambda/\fn_\lambda^n)$ is surjective so that $R^1\varprojlim_n$ vanishes.}
\[ 0\to R^1\varprojlim_n \pi_{i+1}F(R_\lambda/\fn_\lambda^n)  \to 
\pi_i\varprojlim_nF(R_\lambda/\fn_\lambda^n) \to \varprojlim_n \pi_iF(R_\lambda/\fn_\lambda^n) \to 0,\]
the third from (Fin), and the last from (Val)$_N$ since $\Rinf/\fn$ is a valuation ring. This completes the proof.
\end{proof}

\begin{rema}\label{rem2;theo:procdhWebelVanishing}
The proof of Theorem \ref{theo:procdhWebelVanishing} shows
$E_2^{p,q} =0$ for $p>d$ and $q>N$ in \eqref{DescentSS}.
This implies a natural isomorphism
\eq{procdhWebelVanishing-isom}{ \pi_{-N-d}F(X) \simeq H^d_{\procdh}(X,a_{\procdh}\pi_{-N}\Fc).}
\end{rema}

In the rest of this section, we will produce a number of presheaves that satisfy the conditions of Theorem~\ref{theo:procdhWebelVanishing}, see \eqref{KtheoryWeibel},
\eqref{NilhLOWeibel}, \eqref{NilHNWeibel}, \eqref{NilZnsynWeibel}, \eqref{NilTCWeibel}, \eqref{ZnEMWeibel}.
These vanishing results are not new and known from \cite{KST-Weibel} and \cite{EM23}.
Theorem \ref{theo:procdhWebelVanishing} gives a new proof of these results.

\subsection{Algebraic $K$-theory}

Non-connective $K$-theory is viewed as an object of $\PSh(\fSchffpS,\Sp)$ using the equivalence \eqref{Shvprocdhcont}.

\begin{lemm} \label{lemm:KTheory}
Non-connective $K$-theory satisfies (Desc), (Fin), (Rig)$_0$ and (Val)$_0$.
\end{lemm}
\begin{proof}
(Desc) follows from \cite[Th.A]{KST-Weibel}, (Fin) is well-known, (Rig)$_0$ follows from the nil-invariance of negative $K$-theory and (Val)$_0$ follows from \cite[Th. 1.3(iii)]{KM21}. 
\end{proof}

From Theorem \ref{theo:procdhWebelVanishing} and Lemma \ref{lemm:KTheory}, we get
\eq{KtheoryWeibel}{K_i(X)=0\qfor i<-d,}
for a Noetherian scheme $X$ of finite Krull dimension $d$, which is the celebrated theorem \cite[Th.B]{KST-Weibel}.
Moreover, \eqref{procdhWebelVanishing-isom} implies
\[ K_{-d}(X) \simeq H^d_{\procdh}(X,a_{\procdh}\pi_0\Kc) ,\]
where $\Kc$ is the presheaf on $\fSchffpX$ associated to $K$ by the formula
\eqref{eq2;continutity}.
Now, we prove that there is a natural isomorphism
\eq{procdh-cdh}{ H^i_{\procdh}(X,a_{\procdh}\pi_0\Kc) \simeq H^i_{\cdh}(X,\ZZ)
\qfor i\geq 0 }
to get a natural isomorphism $K_{-d}(X) \simeq H^d_{\cdh}(X,\ZZ)$, which is \cite[Cor.D]{KST-Weibel}.
Letting $\cC=\cD(\ZZ)$, we have a pair of adjoint functors (cf. \eqref{adjunction-iota})
\begin{equation*}
\cont: \PSh(\Schft_X,\cC)\begin{smallmatrix}\longleftarrow \\ \longrightarrow\\
\end{smallmatrix} \PSh(\fSchffp_X,\cC):\iota^* 
\end{equation*}
inducing an equivalence (cf. \eqref{Shvprocdhcont})
\[ \cont: \Shvzen(\Schft_X,\cC) \isom \Shvzenc(\fSchffp_X,\cC).\]
So, it is enough to show that there is a natural equivalence in $\PSh(\fSchffp_X,\cC)$:
\[\theta:  a_{\procdh} \pi_0\Kc \isom (a_{\cdh} \pi_0 K)^{\cont}=(a_{\cdh} \ZZ)^{\cont}.\]
To construct a functor $\theta$, it is enough to define
$\pi_0\Kc \to (a_{\cdh} \pi_0 K)^{\cont}$ since $(a_{\cdh} \pi_0 K)^{\cont}\in \Shvzenc(\fSchffp_X,\cC)$. By the above adjunction, it is enough to define
$\iota^*\pi_0\Kc \to a_{\cdh} \pi_0 K$, which follows from the fact that $\iota^*\pi_0\Kc=\pi_0 K$.
To show that $\theta$ is an equivalence, note that both domain and target of $\theta$ are in $\Shvzen(\fSchffp_X,\cC)$, By Theorem~\ref{prop:procdhLocal}, it is enough to show 
$\varinjlim_\lambda\theta(R_\lambda)$ is an equivalence 
for every ind-system $\catinjlim \lambda R_\lambda$ as in the proof of Theorem \ref{theo:procdhWebelVanishing}.
We have
\[ \varinjlim_\lambda (a_{\procdh} \pi_0\Kc)(R_\lambda)=\varinjlim_\lambda (\pi_0\Kc)(R_\lambda) 
 \overset{\eqref{eq2;continutity}}{=} \varinjlim_\lambda \pi_0 \varprojlim_n K(R_\lambda/\fn_\lambda^n) =  \ZZ,
\]
where the last equality follows from Milnor's exact sequence
\[0\to R^1\varprojlim_n K_1(R_\lambda/\fn_\lambda^n)  \to 
\pi_0\varprojlim_n K(R_\lambda/\fn_\lambda^n) \to \varprojlim_n K_0(R_\lambda/\fn_\lambda^n) \to 0,\]
and the nil-invariance (resp. nil-surjectivity) of $K_0$ (resp. $K_1$) and also from the fact that $\varinjlim_\lambda \Rl/\fn_\lambda$ is local. On the other hand, we have
\begin{multline*}
\varinjlim_\lambda (a_{\cdh} \pi_0 K)^{\cont}(R_\lambda) \overset{\eqref{eq2;continutity}}{=}
 \varinjlim_\lambda \varprojlim_n (a_{\cdh}\pi_0 K)(R_\lambda/\fn_\lambda^n)=  
\varinjlim_\lambda (a_{\cdh} \pi_0 K)(R_\lambda/\fn_\lambda)\\
=  (a_{\cdh} \pi_0 K)(\varinjlim_\lambda \Rl/\fn_\lambda)=   \ZZ,
\end{multline*}
where the second equality follows from the nil-invariance of cdh sheaves, the third from the finitarity of $a_{\cdh} \ZZ$ and the last from the fact that $\varinjlim_\lambda \Rl/\fn_\lambda$ is a valuation ring.
These imply the desired equivalence and completes the proof of \eqref{procdh-cdh}.

\subsection{Negative cyclic homology} %
\newcommand{\Fmot}{Fil_{\operatorname{mot}}}
\newcommand{\Znsm}{\ZZ(n)^{\operatorname{sm}}}
\newcommand{\Znem}{\ZZ(n)^{\operatorname{EM}}}
\newcommand{\Zncdh}{\ZZ(n)^{\cdh}}
\newcommand{\Znemc}{\ZZ(n)^{\operatorname{EM,\cont}}}
\newcommand{\Znpcdh}{\ZZ(n)^{\procdh}}
\newcommand{\Lsm}{L^{\operatorname{sm}}}
\newcommand{\Fpcdh}{Fil_{\procdh}}

\newcommand{\HKR}{\operatorname{HKR}}
\newcommand{\FHKR}{\ensuremath{Fil_{\HKR}}}
\newcommand{\BMS}{\operatorname{BMS}}
\newcommand{\FBMS}{\ensuremath{Fil_{\BMS}}}
\newcommand{\gr}{\operatorname{gr}}

\def\Hmot{H_{\operatorname{mot}}}
\begin{defi}\label{def:Nil}
For $\cE \in \PSh(\Sch^{\qcqs}, \cC)$, define
\begin{equation} \label{NilF}
{ \Nil \cE:=\fib(\cE \to a_{\cdh} \cE)\in \PSh(\Sch_\FF^{\qcqs},\cC),}
\end{equation}
where $a_{\cdh}: \PSh(\Sch_\FF^{\qcqs},\cC)\to \Shv_{\cdh}(\Sch_\FF^{\qcqs},\cC)$ is the cdh sheafification functor.
\end{defi}

\begin{rema}\label{rem;def:Nil}
Since $\cE(R)=(a_{\cdh} \cE)(R)$ for any valuation ring $R$, 
$\Nil \cE$ satisfies the following stronger condition than (Val)$_N$ from Theorem \ref{theo:procdhWebelVanishing}:
\begin{itemize}
\item[(Val)] 
For any valuation ring $R$, $F(R)=0$.
\end{itemize}
\end{rema}

Recall that for a $\QQ$-algebra $R$, one defines $\HN(R/\QQ) = \mathrm{HH}(R/\QQ)^{hS^1}$ as the homotopy fixed points of the Hochsdhild homology and extends it to $X\in \Sch^{\qcqs}_\QQ$ using Zariski descent. In \cite[Thm.1.1]{An}, Antieau defines a functorial complete decreasing multiplicative $\ZZ$-indexed filtration%
%
%
\begin{equation} \label{HKR}
\biggl\{\FHKR^n \HN(X/\QQ)\biggr\}_{n\in \ZZ}\;\text{on $\HN(X/\QQ)$}
\end{equation}
and natural equivalences  
\begin{equation*} \label{equa:grFHKR}
\gr_{\FHKR}^n \HN(X/\QQ) \simeq \hLOH n {X/\QQ}[2n].
\end{equation*}
Here, $\hLO {X/\QQ}$ is the Hodge-completed derived de Rham complex equipped with the Hodge filtration $\bigl\{\hLOH n {X/\QQ}\bigl\}_{n\in \NN}$.
We have an equivalence in $ \PSh(\Sch^{\qcqs}_{\QQ},D(\QQ))$:
\[\Nil\FHKR^0\HN(-/\QQ) \stackrel{\sim}{\to} \Nil\HN(-/\QQ) \]
since the cofibre ${\operatorname{cofib}}(\FHKR^0\HN(-/\QQ) \to \HN(-/\QQ))$ is a cdh sheaf on $\Sch^\qcqs_\QQ$ by \cite[Lem.4.5]{EM23}.
Therefore, applying $\Nil$ to \eqref{HKR} produces a complete and exhaustive $\NN$-indexed filtration 
\begin{equation} \label{NilHKR0}
{\biggl\{\FHKR^n \Nil \HN(-/\QQ)\biggl\}_{n\in \NN}\;\text{ on $\Nil\HN(-/\QQ)$}} 
\end{equation}
with identifications
\begin{equation*} \label{NilHKR}
{ \gr^n_{\FHKR} \Nil \HN(-/\QQ) \simeq \Nil\hLOH n {-/\QQ}[2n].}
\end{equation*}

\begin{lemm} \label{lemm:nilhLO}
$\Nil \hLO {-/\QQ}^{\geq n} \in \PSh(\Sch^\qcqs_\QQ,D(\QQ))$ satisfy (Desc), (Fin), (Rig)$_{n} $.
\end{lemm}
\begin{proof}
This follows from \cite[Prop.8.12]{EM23}, see also \cite[Lem.8.7]{KS23} for a proof of the conditions (Desc) and (Fin).
Here, we include a proof of (Rig)$_{n} $. 
Let $A$ be a Noetherian $\QQ$-algebra and $I\subset A$ be a nilpotent ideal.
By \cite[Lem.4.5]{EM23}, $\hLO {-/\QQ}$ is a cdh sheaf so that 
$\Nil \hLO {-/\QQ} = 0$ leading to an equivalence 
\[\Nil \hLO {-/\QQ}^{\geq n} \cong \Nil \LO {-/\QQ}^{<n}[-1].\]
So it suffices to show that the fibre of
$ \Nil \LO {A/\QQ}^{<n}\to \Nil \LO {(A/I)/\QQ}^{<n}$ is supported in cohomological degrees $\leq n-1$. Since $G(A)\simeq G(A/I)$ for any cdh sheaf $G$, it is enough to show that the fibre of
$ \LO {A/\QQ}^{<n}\to \LO {(A/I)/\QQ}^{<n}$ is supported in cohomological degrees $\leq n-1$. 
Since $\LO {B/\QQ}^{<n}$ for any $\QQ$-algebra $B$ is supported in degrees $\leq n-1$, we are reduced to showing the surjectivity of the map
\[ H^{n-1}(\LO{A/\QQ}^{<n}) \to H^{n-1}(\LO{(A/I)/\QQ}^{<n}) .\]
This holds since the map is identified with
$\Omega^{n-1}_{A/\QQ}\to \Omega^{n-1}_{(A/I)/\QQ}$.
\end{proof}

As a corollary of Theorem \ref{theo:procdhWebelVanishing}, Lemma \ref{lemm:nilhLO} and Remark \ref{rem;def:Nil}, 
we conclude that for a Noetherian $\QQ$-scheme $X$ of finite Krull dimension $d$, we have   
\eq{NilhLOWeibel}{H^i(\hLO {X/\QQ}^{\geq n}) =0\qfor i>d+n,}
from which we deduce 
\eq{NilHNWeibel}{\pi_i\Nil\HN(X/\QQ)=0\qfor i<-d}
using the spectral sequence induced by \eqref{NilHKR0}:
\[ E_2^{i,j}= H^{i-j}(\Nil\hLOH {-j} {X/\QQ}) \Rightarrow H^{i+j}\Nil\HN(X/\QQ).\]

\subsection{Integral topological cyclic homology}

For $X \in \Sch^{\qcqs}_{\FF_p}$, write $\TC(X)$ for the integral topological cyclic homology of $X$. By \cite{BMS2}, there exists a functorial complete decreasing $\NN$-indexed filtration
\begin{equation} \label{eq;BMSfilt}
{\biggl\{ \FBMS^n \TC(X)\biggl\}_{n\in \NN} \text{ on } \TC(X)} 
\end{equation}
with associated graded quotients 
\[\gr_{\FBMS}^n \TC(X) \simeq \Znsyn(X)[2n]\]
for a natural object $\Znsyn\in \PSh(\Sch^{\qcqs}_{\FF_p},D(\ZZ_p))$ called the syntomic complex. 

\begin{lemm} \label{lemm:NilZsyn}
$\Nil\Znsyn \in \PSh(\Sch^{\qcqs}_{\FF_p}, D(\ZZ_p))$ satisfy (Desc), (Fin), (Rig)$_{n} $.
\end{lemm}
\begin{proof}
This follows from \cite[Prop.8.12]{EM23}, see also \cite[Lem.8.9]{KS23} for a proof of the conditions (Desc) and (Fin).
Here, we include a proof of (Rig)$_{n} $. 
Let $A$ be a Noetherian $\FF_p$-algebra and $I\subset A$ be a nilpotent ideal.
Since $G(A)\simeq G(A/I)$ for any cdh sheaf $G$, we have
\[\fib(\Nil\Znsyn(A)\to \Nil\Znsyn(A/I))=\fib(\Znsyn(A)\to \Znsyn(A/I)).\]
Hence, (Rig)$_{n} $ follows from \cite[Th.5.2]{AMMN}.
\end{proof}

As a corollary of Theorem \ref{theo:procdhWebelVanishing}, Lemma \ref{lemm:NilZsyn} and Remark \ref{rem;def:Nil}, 
we conclude that for a Noetherian $\FF_p$-scheme $X$ of finite Krull dimension $d$, we have 
\eq{NilZnsynWeibel}{H^i(\Nil\Znsyn(X))=0\qfor i>d+n.}
from which we deduce 
\eq{NilTCWeibel}{\pi_i\Nil\TC(X)=0\qfor i<-d.}
using the spectral sequence induced by \eqref{eq;BMSfilt}:
\[ E_2^{i,j}= H^{i-j}(\Nil\ZZ_p(-j)^{\mathrm{syn}}(X)) \Rightarrow \pi_{-i-j}\Nil\TC(X).\]

\subsection{Motivic cohomology} %
\newcommand{\noe}{\mathrm{noe}}
In what follows, $\FF=\QQ$ or $\FF_p$ and let $\Sm_\FF$ denote the category of  smooth schemes separated of finite type over $\FF$.
Let $\Znsm$ be Voevodsky's $\AA^1$-invariant motivic complex defined in \cite{SVBK} as 
\begin{equation*} \label{equa:ZZtrGm}
\Znsm(X) = 
\underline{C}_*(\ZZ_{tr}(\GG_m^{\wedge n}))(X)[-n] \qfor X\in \Sm_\FF.
\end{equation*}
This is strictly functorial in $X\in \Sm_\FF$. 
In \cite[Cor.2]{Voe02}, this is shown to be quasi-isomorphic scheme-wise to Bloch's cycle complex. 
Recently, Elmanto-Morrow \cite{EM23} introduced a new motivic complex defined for all $X\in \Sch^\qcqs$:
\[ \Znem\in \PSh(\Sch^\qcqs_\FF,D(\ZZ))).\]
They construct $\Znem$ by modifying the cdh sheafification $\Zncdh=a_{\cdh} \Lsm \Znsm$ of the left Kan extension $\Lsm \Znsm$ of $\Znsm$ along $\Sm_\FF \to \Sch^\qcqs_\FF$ by using Hodge-completed derived de Rham complexes in case $\FF=\QQ$ and syntomic complexes in case $\FF=\FF_p$.
The construction is motivated by trace methods in algebraic $K$-theory using the cyclotomic trace map. In \cite[\S9]{KS23}, we gave a different construction of a motivic complex:
\[\Znpcdh:=a_{\procdh} \Lsm \Znsm \in \Shv_{\procdh}(\Sch^\qcqs_\FF,D(\ZZ)),\]
as the sheafification of $\Lsm \Znsm$ with respect to the pro-cdh topology \cite[Def.1.1]{KS23}.
It is shown \cite[Cor.1.11]{KS23} that for any Noetherian $\FF$-scheme $X$ with $\dim(X)<\infty$, there are equivalences
\[ \Znpcdh(X) \simeq \Znem(X) \]
functorial in $X$. We then define a motivic cohomology of such $X$ as
\[   \Hmot^i(X,\ZZ(n)) = H^i( \Znem(X)) =H^i(\Znpcdh(X) ).\]

\begin{lemm} \label{lemm:Znem}
$\Znem \in \PSh(\Sch^{\qcqs}_{\FF}, D(\ZZ))$ satisfy (Desc), (Fin), (Rig)$_{n} $, (Val)$_{n} $.
\end{lemm}
\begin{proof}
This follows from \cite[Prop.8.12]{EM23}. Here, we include a proof.
By \cite[Th.4.10(2), Th.4.24.(2)]{EM23}, there are fiber sequences in $ \PSh(\Sch^\qcqs_\FF,D(\ZZ)))$
\[ \Nil\hLOH n {-/\QQ} \to  \Znem \to \Zncdh\;\text{ if } \FF=\QQ,\]
\[ \Nil\Znsyn \to  \Znem \to \Zncdh\;\text{ if } \FF=\FF_p.\]
By Lemmas \ref{lemm:nilhLO} and \ref{lemm:NilZsyn}, it suffices to show that 
$\Zncdh$ satisfy (Desc), (Fin), (Rig)$_{n} $, (Val)$_{n} $.
(Desc) follows from the fact that the cdh topology is finer than the pro-cdh topology,
(Fin) from that the cdh sheafification of a finitary sheaf is finitary, 
(Rig)$_{n} $ from that any cdh sheaf is invariant with respect to nilpotenet ideals, and
(Val)$_{n} $ from that we have $\Zncdh(R)=\Lsm\Znsm(R)$ for a valuation ring $R$ and the latter is supported in degree$\leq n$ (see \cite[Rem.9.6]{KS23}).
\end{proof}
\def\hK{\hat{K}}

Theorem \ref{theo:procdhWebelVanishing} and Lemma \ref{lemm:Znem} imply that
for a Noetherian $\FF$-scheme $X$,
\eq{ZnEMWeibel}{\Hmot^i(X,\ZZ(n))=0\;\text{  for $i>\dim(X)+n$}}
which is \cite[Th.1.10]{EM23}.



\end{document}